%% file: bitab.tex
\documentclass[10pt]{amsart}
\usepackage{amsfonts,amsmath, amssymb, latexsym, amsthm}
\usepackage{epsf}
\usepackage{graphicx}
\usepackage{bigstrut,multirow,tabularx}

\usepackage[vcentermath]{youngtab} 

\newcommand{\myfig}[3]{\begin{figure}[htbp]
\begin{center}
{\scalebox{#1}{\includegraphics{#2}}}\caption{#3}\label{fig:#2}
\end{center}
\end{figure}}

\input{youngdiag}
\def\la{\lambda}

\def\th{^{th}}

\def\timenow{\count5=\time \divide\count5 by 60 \count6=\time
  \count7=\count5 \multiply\count7 by -60 \advance\count6 by \count7 }

\def\partition{\vdash}

\newcommand{\sumsb}[1]{\sum_{\substack{#1}}}

\DeclareMathOperator{\sgn}{sgn}
\DeclareMathOperator{\per}{per}
\DeclareMathOperator{\pow}{pow}

\newcommand{\bbC}{\mathbb{C}}
\newcommand{\bbN}{\mathbb{N}}
\newcommand{\bbZ}{\mathbb{Z}}

\newcommand{\cxt}{X(\rowseq(\cochg{\gamma}(T)))} 
\newcommand{\cyt}{Y(\rowseq(\cochg{\gamma}(T)))} 

\newcommand{\balpha}{\boldsymbol{\alpha}}
\newcommand{\bbeta}{\boldsymbol{\beta}}
\newcommand{\bdelta}{\boldsymbol{\delta}}
\newcommand{\blambda}{\boldsymbol{\lambda}}
\newcommand{\bnu}{\boldsymbol{\nu}}
\newcommand{\bepsilon}{\boldsymbol{\epsilon}}
\newcommand{\bmu}{\boldsymbol{\mu}}

\newcommand{\latdiag}[1]{L[#1]}             
\newcommand{\hollatdiag}{{L_\gamma}}
\newcommand{\Dalph}{\Delta_{\latdiag{\balpha}}}
\newcommand{\Dgenalph}[1]{\Delta_{\latdiag{#1}}}
\newcommand{\Ialpha}{\mathcal{I}_{\latdiag{\balpha}}}
\newcommand{\Igam}{\mathcal{I}_{\gamma}}

\newcommand{\frch}[1]{\mathcal{F}\,\mathrm{ch}(#1)}
\newcommand{\grch}[1]{\mathrm{ch}(#1)}
\newcommand{\hilb}[1]{\mathcal{H}(#1)}

\newcommand{\bitabsetpn}{\Theta_n'}
\newcommand{\bitabsetn}{\Theta_n}

\newcommand{\colstr}[1]{\mathcal{CS}_{#1}}
\newcommand{\colstrn}[1]{\mathcal{CS}_{n,#1}}
\newcommand{\stdtaball}{SYT}
\newcommand{\stdtabla}{SYT(\lambda)}
\newcommand{\stdtabn}{SYT_n}
\newcommand{\genalph}{\Sigma}
\newcommand{\cochg}[1]{C_{#1}}

\newcommand{\cochgalln}[1]{\mathcal{CO}_{n,#1}}
\DeclareMathOperator{\diagram}{dg}
\DeclareMathOperator{\shape}{sh}
\DeclareMathOperator{\rowseq}{rs}
\DeclareMathOperator{\colseq}{cs}
\DeclareMathOperator{\content}{cont}
\DeclareMathOperator{\std}{std}
\DeclareMathOperator{\conx}{con_X}
\DeclareMathOperator{\cony}{con_Y}

\newcommand{\ol}[1]{\overline{#1}}
\newcommand{\ul}[1]{\underline{#1}}
\newcommand{\ula}{\underline{0}}
\newcommand{\ulb}{\underline{1}}
\newcommand{\ulc}{\underline{2}}
\newcommand{\uld}{\underline{3}}
\newcommand{\ule}{\underline{4}}

\newcommand{\olb}{\bar{1}}
\newcommand{\olc}{\bar{2}}
\newcommand{\old}{\bar{3}}
\newcommand{\ole}{\bar{4}}
\newcommand{\olf}{\bar{5}}

\newcommand{\Wtild}{\tilde{W}}
\newcommand{\Stild}{\tilde{S}}

\newcommand{\hgam}{h_\gamma}

\newcommand{\A}{\mathcal A}
\newcommand{\lexord}[1]{{\mathrm{lex}(#1)}}
\newcommand{\normlexord}{{\mathrm{lex}}}
\newcommand{\bitabord}{<_\mathrm{bitab}}
\newcommand{\bitabordg}{>_\mathrm{bitab}}

\newcommand{\revlexordstr}{<_\mathrm{rlex}}
\newcommand{\detord}{{\mathrm{det}}} 
\newcommand{\perord}{{\mathrm{per}}} 

\newcommand{\CXY}{\bbC[X_n,Y_n]}
\newcommand{\quor}[1]{\CXY_{#1}}
\newcommand{\RI}{\CXY_{\Igam}}
\newcommand{\XYn}{(X_n,Y_n)}
\newcommand{\GS}{{{\mathcal G}_\gamma}\XYn}
\newcommand{\IS}{{{\mathcal I}_\gamma}\XYn}
\newcommand{\JS}{{{\mathcal J}_\gamma}\XYn}
\newcommand{\KS}{{{\mathcal K}_\gamma}\XYn}

\newcommand{\HS}{{{\mathcal H}_\gamma}\XYn}

\newcommand{\ISnxy}{{\mathcal I}_\gamma}
\newcommand{\GSnxy}{{\mathcal G}_\gamma}
\newcommand{\JSnxy}{{\mathcal J}_\gamma}
\newcommand{\KSnxy}{{\mathcal K}_\gamma}

\newcommand{\HSnxy}{{\mathcal H}_\gamma}

\newcommand{\tildlam}{\tilde{\lambda}}

\def\CXY{\C[X_n,Y_n]}
\def\DS{\Delta_{\gamma}(X,Y)}
\def\DSnoxy{\Delta_{\gamma}}

\newcommand{\prt}[1]{\partial_{#1}}

\def\basseqold{\left[\left[1^{m_1},k_1,1^{p_1}\right],
    \left[1^{m_2},k_2,1^{p_2}\right]\right]}

\def\basseq{(m,k,p)}

\def\Rstar{\mathcal{R}^*}
\def\BPgam{\mathcal{B}_\gamma}
\def\hhBPgam{\mathcal{BB}_\gamma}
\def\hhBPgamii{\mathcal{BB}_{\gamma_1,\gamma_2|}}

\def\BPbasis{\mathcal{BP}}
\def\BDbasis{\mathcal{BD}}
\def\EEBbasis{\mathcal{EEB}}

\def\RJ{{\Rstar_{\JSnxy}}}

\def\RH{{\Rstar_{\HSnxy}}}
\def\RG{{\Rstar_{\GSnxy}}}

\def\C{{\hbox{\Ch C}}}
\def\C{{\mathbb{C}}}
\def\Rstar{\C[X_n,Y_n]}

\theoremstyle{plain} \newtheorem{theorem}{Theorem}
\newtheorem{lemma}[theorem]{Lemma}

\newtheorem{corollary}[theorem]{Corollary}

\theoremstyle{definition} 

\newtheorem{example}[theorem]{Example}
\newtheorem{note}[theorem]{Remark}

\begin{document}
\title{ Bitableaux Bases for Garsia-Haiman Modules of Hollow type}

\author{Edward E. Allen}
\address{Department of Mathematics, Wake Forest University, 
  Winston-Salem, NC}
\email{allene@wfu.edu}

\author{Miranda E. Cox}
\address{Section on Biostatistics, Department of Public Health Sciences, 
  Wake Forest University School of Medicine, Winston-Salem, NC}
\email{micox@wfubmc.edu}

\author{Gregory S. Warrington}
\address{Department of Mathematics, Wake Forest University, 
  Winston-Salem, NC}
\email{warrings@wfu.edu}

\date{\today}

\begin{abstract}
  Garsia-Haiman modules $\RI$ are quotient rings in variables
  $X_n=\{x_1, x_2, \ldots, x_n\}$ and $Y_n=\{y_1, y_2, \ldots, y_n\}$
  that generalize the quotient ring $\C[X_n]/{\mathcal I}$, where
  ${\mathcal I}$ is the ideal generated by the elementary symmetric
  polynomials $e_j(X_n)$ for $1 \le j \le n$.  A bitableau basis for
  the Garsia-Haiman modules of hollow type is constructed.
  Applications of this basis to representation theory and other
  related polynomial spaces are considered.
\end{abstract}

\maketitle

\section{Introduction}

Let $X_n = \{x_1,\ldots,x_n\}$ and $Y_n = \{y_1,\ldots,y_n\}$ be sets of
indeterminates.  The main purpose of this paper is to give explicit
combinatorial bases for certain quotients of the ring
\begin{equation}
  \CXY = \bbC[x_1, x_2, \ldots, x_n, y_1, y_2, \ldots, y_n]
\end{equation}
of polynomials in the variables $X_n$ and $Y_n$ with complex
coefficients.  In doing so, we give combinatorial interpretations
for the corresponding Hilbert and Frobenius series.  The ideals in the
aforementioned quotients are defined via determinants as described
below.

Throughout this paper, we will identify any element
$\alpha_i=(\alpha_{i,1},\alpha_{i,2})\in \bbN^2$ with the unit square
in the first quadrant of the plane having $\alpha_i$ as its corner
closest to the origin.  A \emph{lattice diagram}, $\latdiag{\balpha} =
(\alpha_1,\ldots,\alpha_n)$, is a sequence of such unit squares.
Writing $z_j^{\alpha_i}$ for the product
$x_j^{\alpha_{i,1}}y_j^{\alpha_{i,2}}$, to any lattice diagram
$\latdiag{\balpha}$ we associate a determinant
\begin{equation}
  \Dalph = \Dalph(X_n,Y_n) = \det
  \begin{pmatrix}
    z_1^{\alpha_1} & z_2^{\alpha_1}& \cdots & z_n^{\alpha_1}\\
    z_1^{\alpha_2} & z_2^{\alpha_2}& \cdots & z_n^{\alpha_2}\\
    \vdots & \vdots & \ddots & \vdots \\
    z_1^{\alpha_n} & z_2^{\alpha_n}& \cdots & z_n^{\alpha_n}
  \end{pmatrix}.
\end{equation}
Given any polynomial $P(X_n,Y_n)\in\CXY$, there is a
corresponding polynomial of differential operators
\begin{equation}
  P(\partial_X,\partial_Y)  = 
  P(\partial_{x_1},\partial_{x_2},\ldots,\partial_{x_n},\partial_{y_1},
  \partial_{y_2},\ldots,\partial_{y_n}).
\end{equation}
(We write $\partial_{x_i}$ as shorthand for
$\partial/{\partial_{x_i}}$.)  With $\balpha$ as above, define the
ideal
\begin{equation}\label{E:defIS}
  \Ialpha  = \left\{ P(X_n,Y_n)\in \Rstar:
  P(\partial_{X},\partial_{Y}) \Dalph  =  0\right\}
\end{equation}
and write $\quor{\Ialpha}$ for the quotient ring $\CXY/\Ialpha$.

These quotients, known as \textit{Garsia-Haiman modules}, were introduced
by A. Garsia and M. Haiman in \cite{GH}.  A good overview of the
subject can be found in \cite{cdm}.  A. Garsia and M.  Haiman
introduced modules of this type to study the $q,t$-Kostka
coefficients.  This paper will henceforth concern itself
only with Garsia-Haiman modules arising from \emph{hollow lattice
  diagrams}.  Roughly, a hollow lattice diagram is a subset of a hook
shape obtained by removing a (perhaps trivial) contiguous region of
cells from each of the arm and leg of the hook (see
Figure~\ref{fig:hollow}).  More precisely we parametrize hollow
lattice diagrams $\hollatdiag$ by sequences of three pairs $\gamma =
\basseq$, with $m = (m_1,m_2)\in\bbZ^2_{\geq 1}, k =
(k_1,k_2)\in\bbZ_{\geq 1}^2$ and $p = (p_1,p_2)\in \bbN^2$, by setting
$\hollatdiag = \latdiag{\balpha}$ where
\begin{align*}
\balpha &= \left((0,m_2+k_2+p_2-1),(0,m_2+k_2+p_2-2),\ldots,(0,m_2+k_2),\right.\\
    &\phantom{=\ }(0,m_2+k_2-1), (0,m_2-1),\ldots,(0,1),(0,0),(1,0),(2,0),\\ 
    &\phantom{=\ }\left. \ldots,(m_1-1,0),(m_1+k_1-1,0),
      \ldots,(m_1+k_1+p_1-1,0)\right). 
\end{align*}
(We also allow $k_i = 0$ if $p_i = 0$.)  Unless otherwise noted, the
number of cells in $\hollatdiag$ (namely, $m_1 + p_1 + m_2 + p_2 + 1$)
will be denoted by $n$.

\myfig{.35}{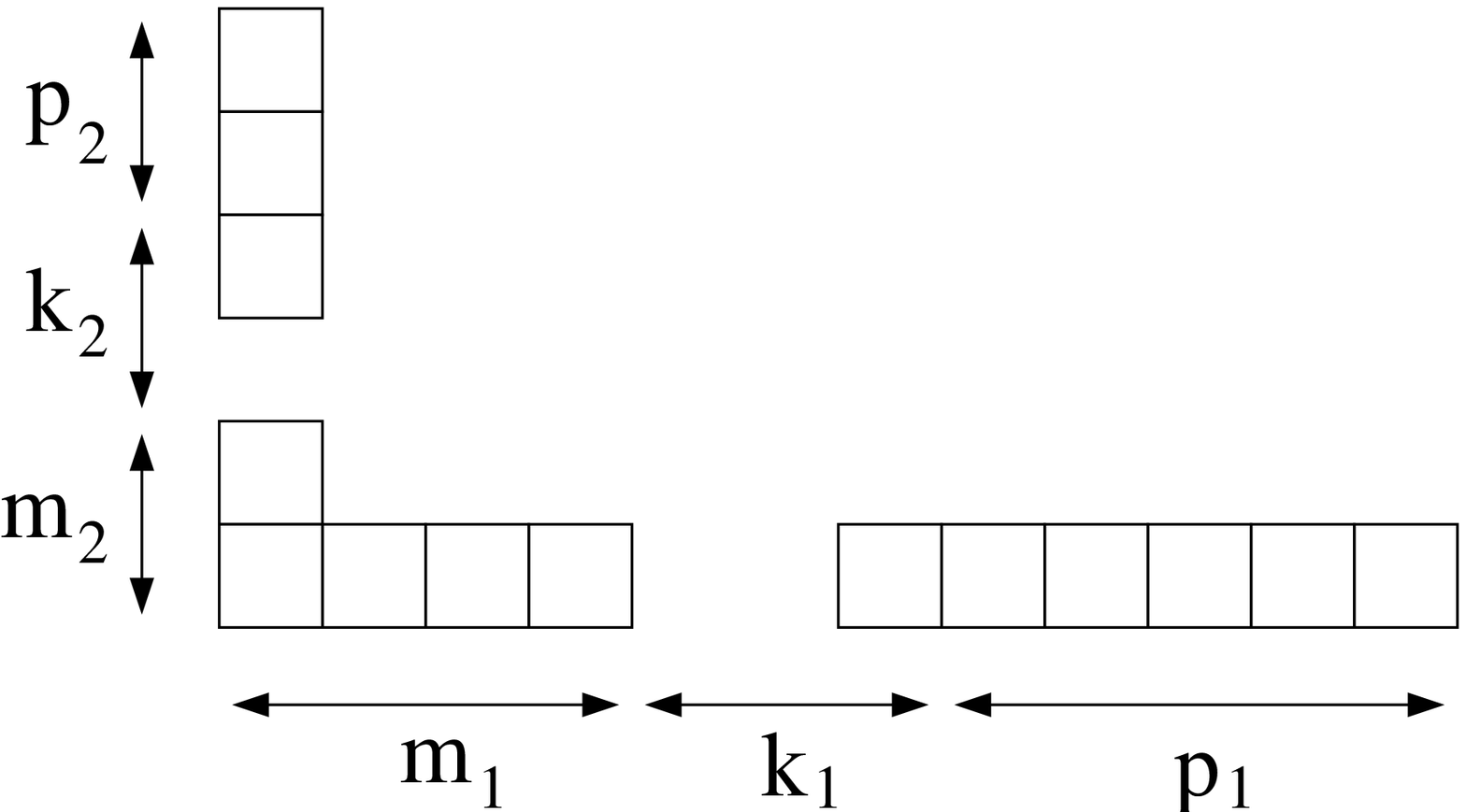}{Hollow lattice diagram $\hollatdiag$ with $m =
  (4,2)$, $k = (3,2)$, and $p = (5,2)$.}

Abusing notation slightly, we write $\Igam$ for
$\mathcal{I}_{\hollatdiag}$ and $\DSnoxy$ for $\Delta_\hollatdiag$.
Our goal is to consider the combinatorics of the hollow Garsia-Haiman
space $\RI$.  As suggested by the previous terminology, the rings
$\RI$ carry $S_n$-representations: The symmetric group $S_n$ has a
natural \emph{diagonal action} on $\CXY$ given by
\begin{equation}
  \sigma P(x_1,\ldots, x_n, y_1,\ldots, y_n) = 
  P(x_{\sigma(1)},\ldots, x_{\sigma(n)}, y_{\sigma(1)},\ldots, y_{\sigma(n)}).
\end{equation}
This action passes through to an action on each $\RI$.  

Let $R$ be any $S_n$-module realized as a polynomial ring over $X_n$
and $Y_n$ (such as $\RI$) and let $R_{r,s}$ denote the subspace of $R$
containing elements of total degree $r$ in $X_n$ and total degree $s$
in $Y_n$.  We can decompose each $R_{r,s}$ as $R_{r,s} = \oplus_{r,s}
c_{r,s} S^\lambda$ where each $S^\lambda$ is an irreducible
$S_n$-module (i.e., Specht module).  Denote the Schur functions by
$s_\lambda$.  The \emph{(bi-graded) character}, \emph{Frobenius
  series} and \emph{Hilbert series} are then, respectively, given by
\begin{equation}
\begin{aligned}
   \grch{R} &= 
  \sum_{r,s} \left(\sum_{\lambda\vdash n} c_{r,s} \chi^\lambda\right) \ t^rq^s,\\
   \frch{R} &= 
  \sum_{r,s} \left(\sum_{\lambda\vdash n} c_{r,s} s_\lambda\right) \ t^rq^s,\\
   \hilb{R} &= 
  \sum_{r,s} \dim(R_{r,s})\ t^rq^s.
\end{aligned}
\end{equation}
Here $\chi^\lambda$ denotes the character of $S^\lambda$ and
$\lambda\vdash n$ signifies that $\lambda$ is a partition of $n$.  The
Frobenius series is the image of the graded character under the
Frobenius map which sends $\chi^\lambda$ to $s_\lambda$.  Note that
the Hilbert series can be recovered from the Frobenius series by
formally replacing each $s_\lambda$ by the dimension of $S^\lambda$.

By constructing an appropriate basis we will prove the following
theorem.  (The definition of a standard tableau will be given in
Section 2; the cocharge statistics $|\cxt|$ and $|\cyt|$ are defined
in Section \ref{sec:HR}.)

\begin{theorem}\label{T:maintheorem}
  Let $\RI$ denote a hollow Garsia-Haiman module parametrized by
  $\gamma = (m,p,k)$. The graded character $\grch{\RI}$ of $\RI$ is given by
  \begin{equation}
    \left[
      \begin{matrix}
        p_1+k_1\\ p_1+1\cr
      \end{matrix}
    \right]_t\ \left[
      \begin{matrix}
        p_2+k_2\\ p_2+1\cr
      \end{matrix}
    \right]_q\ \sum_{\la\partition n} \chi^\la 
    \sum_{T\in \stdtabla} t^{|\cxt|}\ q^{|\cyt|},
  \end{equation}
  where $\stdtabla$ denotes the collection of standard tableau of
  shape $\la$.
\end{theorem}


The modified Macdonald polynomials are a family of symmetric functions
over the field of Laurent polynomials in two variables that specialize
to many important classical symmetric functions.  Parametrized by
partitions, they are given plethystically by
\begin{equation} {\tilde H}_\mu(z;q,t) = J_\mu\left[ \frac{Z}
  {1-t^{-1}} ;q,t^{-1}\right]t^{n(\mu)} = \sum_{\la}\ {\tilde
    K}_{\la,\mu}(q,t) \ s_\la(z).
\end{equation}
An argument analogous to the one used to prove Theorem 2.1 in
\cite{BGT} can be used to deduce the following corollary from
Theorem~\ref{T:maintheorem}.

\begin{corollary}\label{T:Garsiathree}
  The graded Frobenius characteristic of the hollow Garsia-Haiman
  module $\RI$ is given by the polynomial
  \begin{equation} {\mathcal F}\ \grch{\RI}  =  \Xi_{\la}(t)\
    \Xi_{\theta}(q)\ {\tilde H} _{(m_2+p_2+1,1^{m_1+p_1})}(z;q,t)
  \end{equation}
  where $\la=(k_1^{p_1+1})$, $\theta=(k_2^{p_2+1})$ and $\Xi_\nu(t)$
  gives the Hilbert polynomial of the graded vector space of skew
  Schur functions $s_{\nu/\mu}$ as $\mu$ varies in $\nu$.
\end{corollary}

To prove Theorem~\ref{T:maintheorem}, we will define a sequence of
ideals
\begin{equation}\label{eq:ideals}
  \GS \subset \HS \subset \JS \subset \KS \subset \IS.
\end{equation}
For each ideal $\mathcal{E} = \mathcal{E}(X_n,Y_n)$ in
\eqref{eq:ideals}, we will define appropriate generators for
$\mathcal{E}$, construct a base for the corresponding quotient space
$\CXY_{\mathcal{E}} = \CXY/\mathcal{E}$, and compute the corresponding
Hilbert series.  To complete the proof, we will use a correspondence
between our basis elements for $\RI$ and irreducible characters of
$S_n$.  Sections 2, 3 and 4 will introduce the necessary background
and notation on tableaux, cocharge tableaux and symmetric polynomials,
respectively.  The following four sections consider the situations
corresponding to $\GSnxy$, $\HSnxy$, both $\JSnxy$ and $\KSnxy$, and
both $\JSnxy$ and $\ISnxy$, respectively.

We note here that Garsia-Haiman modules corresponding to specific
classes of lattice diagrams have been studied elsewhere.
\emph{Periodic Garsia-Haiman modules} were considered by the first
author \cite{Allena} and (in one variable) by H. Morita and H.-F.
Yamada \cite{Morita}, R. Stanley \cite{Stanley} and J. Stembridge
\cite{Stembridge}.  \emph{Dense Garsia-Haiman modules} were
investigated by the first author in \cite{Allenb}.  The Garsia-Haiman
modules corresponding to the degenerate hollow lattice diagrams of
one-row skew shapes or one-column skew shapes were studied by F.
Bergeron, A. Garsia and G. Tesler in \cite{BGT}.

Finally, it should be noted that a conjecture has been announced by M.
Haiman, J. Haglund, N. Loehr, J.  Remmel, and A. Ulyanov (see
\cite{Haiman}) for a combinatorial formula for the character of the
coinvariants of the symmetric group.

\begin{note}
  There is a wide variety of indexing and notational conventions among
  papers in this field.  Most obviously, when the primary lattice
  diagrams under consideration are partitions, the correspondence
  between $\bbN^2$ and first quadrant lattice points usually has the
  first index giving the $y$-coordinate.  Also note that we here write
  $\basseq$ for the parametrizing tuple $\basseqold$ of \cite{Allenb}.
\end{note}

\section{Tableaux}

A partition $\mu = (\mu_1,\mu_2,\ldots,\mu_j,\ldots)$ is a (possibly
infinite) sequence of weakly decreasing integers with $j \geq 1$
nonzero terms.  We will not distinguish between partitions with the
same collection of nonzero terms.  Write $|\mu| = \mu_1 + \mu_2 +
\cdots + \mu_j$ for the sum of the parts.  If $|\mu| = n$, then we say
that $\mu$ is a \emph{partition} of $n$ and write $\mu\vdash n$.  The
length, $j$, is denoted $\ell(\mu)$.  If $i$ appears $m_i$ times in
$\mu$ for each $i$, then the tuple $(1^{m_1},2^{m_2},\ldots)$ is
called the \emph{type} of $\mu$.  $\mu^t$ will denote the
\emph{transpose} of $\mu$.

Let $\mu,\lambda \vdash n$.  We use ``$\leq_{\normlexord}$'' on
partitions to denote the lexicographic order.  A (French-style)
\emph{Ferrers diagram of shape $\mu$} is a collection of
left-justified unit squares (``cells'') in the first quadrant with
$\mu_i$ squares in the $i$th row from the bottom.  The \emph{shape} of
a Ferrers diagram $D$, $\shape(D)$, is the partition obtained by
listing the row lengths of $D$.  The notation $\diagram(\mu)$ will be
used to denote the canonical Ferrers diagram of shape $\mu$.

A \emph{$\genalph$-filling} $f$ is a map $f: \diagram(\mu)
\longrightarrow \genalph$ from the cells of a Ferrers diagram to some
totally ordered alphabet $\genalph$.  We will consider fillings with
three different alphabets:
\begin{enumerate}
\item $\A' = \{(a,b): a,b\in \bbN\}$ ordered by $(a_1,b_1) <_{\A'}
  (a_2,b_2)$ whenever
  \begin{enumerate}
  \item $a_1 - b_1 < a_2 - b_2$; or
  \item $a_1 - b_1 = a_2 - b_2$ and $a_1 < a_2$.
  \end{enumerate}
  Geometrically, the order $<_{\A'}$ can be visualized as listing the
  points in the first quadrant by reading down lines $y = x + c$ from
  left to right with successively smaller values of $c$.
\item $\A = \{(a,b)\in \A': a = 0 \text{ or } b = 0\}$ with the order
  $<_{\A}$ induced by $<_{\A'}$.  Note that the elements of $\A$ index
  cells that can appear in a hollow lattice diagram.  For brevity in
  formulas, we sometimes write $\ul{a}$ for $(a,0)$ and $\ol{b}$ for
  $(0,b)$.  The notation is meant to evoke positive and negative
  numbers, respectively, as this interpretation of the elements of
  $\A$ is consistent with $<_{\A}$.
\item $\bbN$ ordered by $0 < 1 < 2 < 3 < \cdots$.
\end{enumerate}

The picture (or pair $(f,\diagram(\mu))$) obtained by placing elements
of $\genalph$ in the cells of a Ferrers diagram of shape $\mu$
according to $f$ is a \emph{$\genalph$-filled diagram}.  When
$\genalph$ is clear (or unimportant), we simply refer to \emph{filled
  diagrams.}  For a filled diagram $U = (f,\diagram(\mu))$, we use the
shorthand $g(U)$ for the new filled diagram $(g\circ f,\diagram(\mu))$.

A filled diagram is \emph{injective} if the map $f$ is an injective
map.  When $\genalph = \bbZ_{\geq 1}\subset \bbN$, we refer to a
$\Sigma$-filled diagram as a \emph{tableau}.  A filled diagram of
shape $\mu$ is said to be \emph{column strict} if the entries increase
weakly from left to right in each row and increase strictly in each
column from bottom to top.  We will denote the collection of
column-strict tableaux for a given alphabet $\genalph$ by
$\colstr{\genalph}$ (and use $\colstrn{\genalph}$ if we want to
specify the number of boxes).  An injective column-strict tableau with
distinct entries $\{1,\ldots,n\}$ for some $n$ is often referred to as
a \emph{standard} tableau.  Let $\stdtaball$, $\stdtabla$, and
$\stdtabn$ denote the collections of standard tableaux, standard
tableaux of shape $\lambda$, and standard tableaux with $n$ cells,
respectively.  In the context of filled diagrams, $T$ will be reserved
for a standard tableau; $V$ for a column-strict tableau.  Finally, for
any filled diagram $U$, we define $U^t$ to be the tableau obtained by
reflecting $U$ along the line $y=x$.  Figure~\ref{fig:tabex}
illustrates an $\A$-filled diagram $U$ of shape $\shape(U) = (3,2)$
along with its transpose $U^t$.  Note that in this case, $U^t\in
\colstr{\A}$.

\myfig{.6}{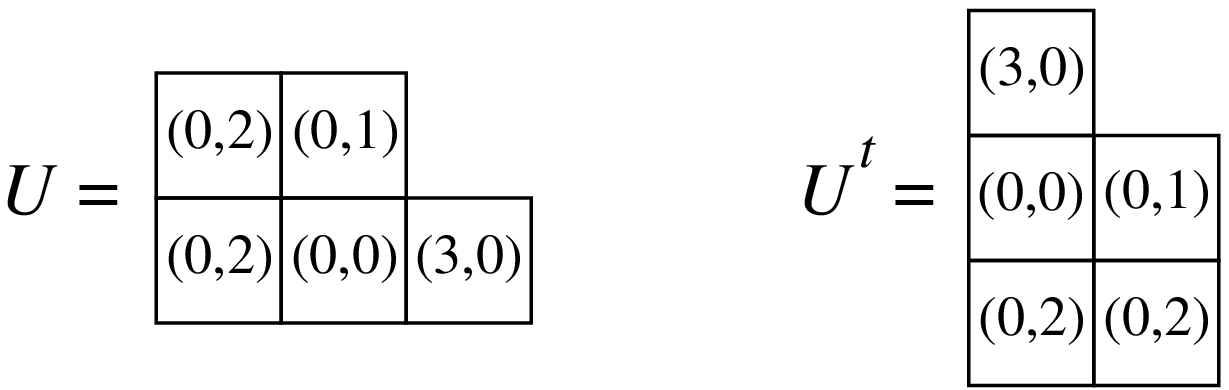}{An $\A$-filled diagram and its transpose.}

Let $I$ denote an injective tableau of shape $\mu=(\mu_1, \ldots,
\mu_j)$, $R_i$ ($1 \le i \le j$) denote the collection of integers in
the $i\th$ row of $I$ and $D_i$ ($1 \le i \le \mu_1$) denote the
collection of integers in the $i\th$ column of $I$. Set
\begin{equation}\label{E:defR}
  R(I)  =  S_{R_1} \times S_{R_2} \times \cdots \times S_{R_j} 
\end{equation}
and
\begin{equation}\label{E:defD}
  D(I)  =  S_{D_1} \times S_{D_2} 
  \times \cdots
  \times S_{D_{\mu_1}}, 
\end{equation}
where $S_{R_i}$ and $S_{D_i}$ denote the symmetric group on the
collections of elements $R_i$ and $D_i$, respectively. Define, in the
group algebra $\C[S_n]$,
\begin{equation}
  P(I) = \sum_{\sigma \in R(I)} \sigma \quad \text{ and }\quad 
  N(I) = \sum_{\sigma \in D(I)} \sgn(\sigma)\ \sigma. 
\end{equation}

A \emph{bitableau} is a pair $(S,U)$ of filled diagrams of the same
shape where $S$ is $\bbN$-filled and $U$ is $\A'$-filled.  A
\emph{standard bitableau} satisfies the additional stipulations that
$S$ is a standard tableau and $U\in\colstr{\A'}$.  The set of all
standard bitableaux $(S,U)$ on $n$ boxes will be denoted $\bitabsetn$
if we restrict to $U\in\colstrn{\A}$ and $\bitabsetpn$ if we do not
place this restriction.

Given any $\genalph$-filled diagram $U = (f,\diagram(\mu))$, define
the \emph{standardization} $\std(U) = (\xi_f,\diagram(\mu))$ by
setting $\xi_f$ to be the unique standard filling such that: For cells
$c,d\in\diagram(\mu)$,
\begin{enumerate}
\item $f(c) \leq_{\genalph} f(d)$ implies $\xi_f(c) < \xi_f(d)$.
\item If $f(c) = f(d)$ and either
  \begin{enumerate}
  \item $c$ is north of $d$, or
  \item $c$ is in the same row as $d$ but west,
  \end{enumerate} 
  then $\xi_f(c) < \xi_f(d)$.
\end{enumerate}
For $U$ a $\genalph$-filled diagram and $T$ a standard tableau, we
denote the entry in $U$ in the cell corresponding to $i$ in $T$ by
$u^T_i$.  Two special cases that arise frequently are when $T =
\std(U)$ and when $(T,U)$ is a bitableau.  In the former case, we
abbreviate $u^T_i$ by $u_i$.

Let $S$ be an injective $\bbN$-filled tableau.  Again writing
$z_j^{u^S_i}$ for $x_j^{u^S_{i,1}}y_j^{u^S_{i,2}}$, for a bitableau
$(S,U)$ we set the {\it bideterminant} $[S,U]_{\det}$ to be
\begin{equation}
  [S,U]_{\det} = N(S)\ z_1^{u^S_1} z_2^{u^S_2} \cdots z_n^{u^S_n}
  = \sum_{\sigma \in D(S)} \sgn(\sigma) 
  z_{\sigma(1)}^{u^S_1} z_{\sigma(2)}^{u^S_2} \cdots z_{\sigma(n)}^{u^S_n}.
\end{equation}
Similarly, the corresponding {\it bipermanent} $[S,U]_{\per}$ is given
by
\begin{equation}
  [S,U]_{\per} = P(S)\ z_1^{u^S_1} z_2^{u^S_2} \cdots z_n^{u^S_n}
  = \sum_{\sigma \in R(S)} 
  z_{\sigma(1)}^{u^S_1} z_{\sigma(2)}^{u^S_2} \cdots z_{\sigma(n)}^{u^S_n}.
\end{equation}
The following theorem is a special case of \cite[Theorem 8]{Rotab}
(also cf.  \cite{Allenb,Rotaa}).  We suggest the reader work out some
examples from the case $n=3$ by hand.
\begin{lemma}\label{T:lemfour}
  The collections 
  \begin{align}\label{E:biPerms}
    \BPbasis &= \left\{[T,V]_{\per}: (T,V)\in\bitabsetpn \right\} \text{ and }\\
    \BDbasis &= \left\{[T,V]_{\det}: (T,V)\in\bitabsetpn \right\}
  \end{align}
  are infinite bases for $\CXY$.
\end{lemma}
For a $\genalph$-filled diagram $U$, the {\it row sequence}
$\rowseq(U)$ is the sequence obtained by listing the entries of $U$ in
each row from left to right, starting with the bottom row.  The
\emph{column sequence} $\colseq(U)$ is found by listing the entries of
$U$ from bottom to top in each column, starting with the leftmost
column.  Finally, the \emph{content} $\kappa(U)$ is a rearrangement of
the row sequence $\rowseq(U)$ of $U$ into nondecreasing order with
respect to $<_{\genalph}$.

\begin{example}\label{ex:rkc}
  For $U$ as in Figure~\ref{fig:tabex}, we have 
  \begin{align*}
    \rowseq(U) &= \left((0,2),(0,0),(3,0),(0,2),(0,1)\right)
    = \left(\ol{2},\ul{0},\ul{3},\ol{2},\ol{1}\right),\\
    \colseq(U) &= \left((0,2),(0,2),(0,0),(0,1),(3,0)\right)
    = \left(\ol{2},\ol{2},\ul{0},\ol{1},\ul{3}\right), \text{ and }\\
    \kappa(U)  &= \left((0,2),(0,2),(0,1),(0,0),(3,0)\right)
    = \left(\ol{2},\ol{2},\ol{1},\ul{0},\ul{3}\right). 
  \end{align*}
\end{example}

There are two orderings of bitableaux that are particularly important
when considering elements of $\BDbasis$ or $\BPbasis$. Let
$>_{\lexord{\A'}}$ denote the lexicographic order with respect to
$>_{\A'}$.  For tableau, $S_1$, $U_1$, $S_2$ and $U_2$, we will say that
\begin{equation}
  (S_1,U_1) <_{\detord} (S_2,U_2)
\end{equation}
whenever
\begin{enumerate}
\item $\shape(S_1^t) <_{\normlexord} \shape(S_2^t)$;
\item If $\shape(S_1) = \shape(S_2)$ then 
  $\kappa(U_1) >_{\lexord{\A'}} \kappa(U_2)$;
\item If $\shape(S_1) = \shape(S_2)$ and $\kappa(U_1) = \kappa(U_2)$ then
\begin{equation}
  \colseq(S_1) \colseq(U_1) >_{\lexord{\A'}} \colseq(S_2) \colseq(U_2),
\end{equation} where $\colseq(S_i) \colseq(U_i)$
is the concatenation of $\colseq(S_i)$ and $\colseq(U_i)$ for $i=1,2$.
\end{enumerate}
\begin{example}
  Let
\begin{equation*}
  (S_1,U_1) = \left(\young(57,348,126)\,, 
    \young(\ula\ulb,\olb\ula\ulb,\olc\olb\olb)\right)\text{ and }
  (S_2,U_2) = \left(\young(57,348,126)\,, 
    \young(\ulb\ulb,\olb\ula\ulb,\olc\olc\olb)\right).
\end{equation*}
Certainly $\shape(S_1^t) = \shape(S_2^t)$.  However
\begin{equation}
  \kappa(U_1) = (\ol{2},\ol{1},\ol{1},\ol{1},\ul{0},\ul{0},\ul{1},\ul{1}) 
  >_{\lexord{\A}} 
  \kappa(U_2) = (\ol{2},\ol{2},\ol{1},\ol{1},\ul{0},\ul{1},\ul{1},\ul{1}).
\end{equation}
So $(S_1,U_1) <_{\det} (S_2,U_2)$.
\end{example}

Similarly, we will say that
\begin{equation}
  (S_1,U_1) <_{\perord} (S_2,U_2) 
\end{equation}
whenever
\begin{enumerate}
\item $\shape(S_1) <_{\lexord{\A}} \shape(S_2)$;
\item If $\shape(S_1) = \shape(S_2)$ then $\kappa(U_1) 
  <_{\lexord{\A'}} \kappa(U_2)$;
\item If $\shape(S_1) = \shape(S_2)$ and $\kappa(U_1) = \kappa(U_2)$ then
  \begin{equation}
    \rowseq(S_1) \rowseq(U_1) >_{\lexord{\A'}} \rowseq(S_2) \rowseq(U_2),
  \end{equation}
  where $\rowseq(S_i) \rowseq(U_i)$ is the concatenation of
  $\rowseq(S_i)$ and $\rowseq(U_i)$ for $i = 1,2$.
\end{enumerate}

\begin{theorem}[\cite{Allenb,Rotaa,Rotab}]\label{T:Rotaexp}
  Let $[S,U]_{\det}$ be a bitableau on $n$ boxes with $U$ $\A'$-filled
  such that either $S$ is not standard or $U$ is not column-strict.
  Then we can write
  \begin{equation}
    [S,U]_{\det} = \sum_i d_i [T_i,V_i]_{\det}
  \end{equation}
  where, for each $i$, it is true that $d_i\in\bbZ$, $(T_i,V_i)\in\bitabsetpn$,
  \begin{equation*}
    (T_i,V_i) >_{\detord} (S,U), 
  \end{equation*}
  $\kappa(T_i) = \kappa(S)$ and $\kappa(V_i) = \kappa(U)$.  The above
  statements hold, \emph{mutatis mutandi}, for bipermanents and the
  order $>_{\per}$.
\end{theorem}

\begin{example}
  \begin{align*}
    \left[\young(3,21)\,,\young(\ulc,\olb\ulb)\right]_{\det} &= 
    (\varepsilon - (2,3))x_1y_2x_3^2 = x_1y_2x_3^2 - x_1y_3x_2^2\\
    &= \left[\young(3,12)\,,\young(\ulc,\olb\ulb)\right]_{\det} - 
    \left[\young(2,13)\,,\young(\ulc,\olb\ulb)\right]_{\det} -
    \left[\young(3,2,1)\,,\young(\ulc,\ulb,\olb)\right]_{\det}.
  \end{align*}
\end{example}

\section{Cocharge Tableaux}

Given a hollow lattice diagram $\hollatdiag$, to any standard tableau
$T = (f,\diagram(\mu))$, we define the \emph{cocharge diagram},
$\cochg{\gamma}(T) = (\hgam\circ\pi\circ f,\diagram(\mu))$.  Here,
$\pi$ is the usual ``cocharge'' map defined recursively by

\begin{equation}
  \pi(i) = 
  \begin{cases}
    0, & \text{ if } i = 1,\\
    \pi(i-1), & \text{ if $i>1$ occurs weakly southeast of $i-1$ in T,}\\
    \pi(i-1) + 1, & \text{ if $i>1$ occurs weakly northwest of $i-1$ in T}.\\
  \end{cases}
\end{equation}
The map $\pi$ is well-defined on standard diagrams.  We
then define $\hgam$ as the unique order- and cover-preserving map from
$\mathbb{N}$ to $\A$ that sends $f^{-1}(m_2+p_2+1)$ to $(0,0)$.

\myfig{.55}{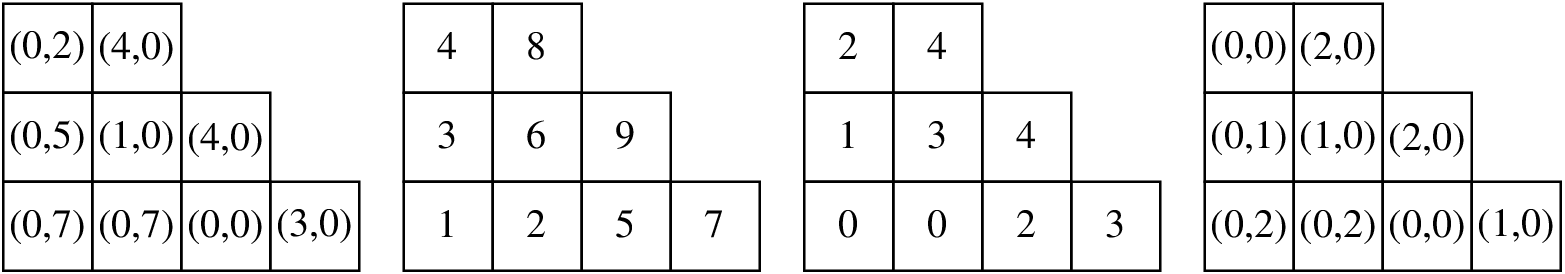}{A diagram $U$, its standardization $\std(U)$,
  $\pi(\std(U))$ and $\cochg{\gamma}(\std(U))$.}

Figure~\ref{fig:cochgII} gives an example cocharge diagram that
corresponds to any $\gamma=(m,k,p)$ describing a lattice diagram with
$9$ boxes such that $m_2 + p_2 + 1$ is equal to $4$ or $5$.

\begin{lemma}\label{lem:decomposition}
  Fix $n$ and $\gamma$.  Define $\cochgalln{\gamma} =
  \{\cochg{\gamma}(T):T\in\stdtabn\}$.  There is a bijection between
  elements $U\in \colstrn{\A}$ with $u_{m_2 + p_2 + 1} = \ula$ and
  pairs $(C,\balpha)$ with $C\in\cochgalln{\gamma}$ and
  $\balpha\in\A^n$ such that $\alpha_i\leq_{\A}\alpha_{i+1}$ for
  $1\leq i\leq n-1$ and $\alpha_{m_2+p_2+1} = (0,0)$.
\end{lemma}
For an $\A$-filled diagram as in Figure~\ref{fig:cochgII}, the sequence
$\balpha$ of Lemma~\ref{lem:decomposition} is
\begin{equation*}
  \balpha = ((0,5),(0,5),(0,4),(0,2),(0,0),(0,0),(2,0),(2,0),(2,0)).
\end{equation*}

\begin{proof}
  Let $T = \std(U)$ and $C = \cochg{\gamma}(T)$.  Map $U$ to
  $(C,\balpha)$ with $\alpha_i = u_i - c^T_i$ for each $1\leq i\leq
  n$.  By the definition of $\hgam$, $c^T_{m_2+p_2+1} = \ula$.
  Combined with our requirement for $u_{m_2+p_2+1}$, it follows that
  $\alpha_{m_2+p_2+1} = \ula$.  By the definitions of $\hgam$, $\pi$
  and standardization, the sequences $c^T_1,c^T_2,\ldots$ and
  $u_1,u_2,\ldots$ are both weakly increasing.  That the $\alpha_i$
  are weakly increasing then follows from the additional fact that $U$
  is column-strict.
\end{proof}

\section{Some Operations by Symmetric Polynomials}

We now review some important definitions and results with respect to
symmetric polynomials. A standard reference for this material is
\cite{Macdonald}.  We introduce the following families of symmetric
polynomials.  We use the convention that tuples of elements of $\A'$
or $\A$ are written in boldface.  In these definitions, let $\lambda =
(\lambda_1,\ldots,\lambda_j,\ldots,\lambda_n)$ be a partition
with $j\leq n$ indexing the last nonzero part.
\begin{enumerate}
\item Define the \emph{monomial symmetric function} as
  \begin{equation}\label{E:monomials}
    m_\la(X_n) = \sum_{\nu=(\nu_1, \nu_2, \ldots, \nu_n)} x_1^{\nu_1}\ x_2^{\nu_2}
    \cdots x_n^{\nu_n}, 
  \end{equation}
  where the sum is over all distinct permutations $\nu$ of
  $\la$.
\item For a sequence $\bbeta = (\beta_1,\ldots,\beta_n)\in
  (\A')^n$, the \emph{MacMahon monomial symmetric function}
  \begin{equation}
    m_{\bbeta}(X_n,Y_n)
    =  \sum_{\bdelta=(\delta_1,\delta_2,\ldots,\delta_n)} 
    z_1^{\delta_1} z_2^{\delta_2} \cdots z_n^{\delta_n},
  \end{equation}
  where the sum is over all distinct permutations $\bdelta$ of $\bbeta$.
\item For a positive integer $r$, the \emph{elementary symmetric
    function}
  \begin{equation}
    e_r (X_n) = \sum_{1\leq i_1 < i_2 < \cdots < i_r \leq n} x_{i_1} x_{i_2} \cdots x_{i_r}.
  \end{equation}
  Set $e_0 = 1$ and $e_\lambda = e_{\lambda_1}\cdots e_{\lambda_j}$.
\item For a positive integer $r$, the \emph{complete (or homogeneous)
    symmetric function}
  \begin{equation}
    h_r(X_n) = \sum_{|\la|=r} m_\la(X_n).
  \end{equation}
  Set $h_0 = 1$ and $h_\lambda = h_{\lambda_1}\cdots h_{\lambda_j}$
  (we extend this definition in the obvious way to the case where
  $\la$ is a $j$-tuple of nonnegative integers; i.e., not necessarily
  nonincreasing).
\end{enumerate}

\begin{lemma}[\cite{Allena}; Theorem 5.4, Corollary 5.5]\label{T:Cortwo}
  Suppose $g: \colstrn{\A'} \rightarrow \colstrn{\A'}$ such that for
  all $W\in\colstr{\A'}$, $W$ and $g(W)$ have the same
  standardization.  Fix $T\in\stdtabn$, $V\in\colstr{\A'}$ and write
  $U = g(V)$.  For each $1\leq i\leq n$, set $\beta_i = v_i - u_i$.
  Write $\bbeta = (\beta_1,\ldots,\beta_n)$.  Suppose that
  $\bbeta\in\A^n$ and $\beta_i\leq_{\A}\beta_{i+1}$ for $1\leq i\leq
  n-1$.  Then
  \begin{multline}\label{eq:mbeta}
    m_{\bbeta}(X_n,Y_n)[T,U]_{\per} = c_{T,V}[T,V]_{\per} + \\
    \sumsb{(S,W) >_{\per} (T,V)\\ S \in \stdtabn,\, W\in\colstrn{\A}}
    c_{S,W}[S,W]_{\per} + 
    \sum_{\Stild \in \stdtabn} c_{\Stild,\Wtild}[\Stild,\Wtild]_{\per},
  \end{multline}
  where $c_{T,V}\ne 0$ and $\Wtild$ is $\A'$-filled with at least one
  entry not in $\A$.
\end{lemma}

\begin{example}
  For $\bbeta = (\olc,\ula,\uld)$, $U = \young(\ulb,\ula\ula)$, $T
  = \young(3,12)$ and $V = \young(\ule,\olc\ula)$, we have
  \begin{equation}
    m_{\bbeta}(X_3,Y_3)\left[T,U\right]_{\per} = 
    2[T,V]_{\per} + 
    2\left[T,\young(\ulb,\olc\uld)\right]_{\per} + 
    2\left[T,\young(a,\ula\uld)\right]_{\per},
  \end{equation}
  where $a = (1,2)$.  We see that the second bipermanent is, in fact,
  larger than $(T,V)$ in the order $>_{\per}$ as the content
  $(\olc,\ulb,\uld)$ is greater than that of $V$ in the lexicographic
  order with respect to $\A$.
\end{example}

We also have the following Lemma (cf. \cite{Allena}, Theorem 5.2).
\begin{lemma}\label{T:expands}
  For $\bbeta\in(\A')^n$ and a lattice diagram $\latdiag{\balpha}$,
  \begin{equation}
    m_{\bbeta}(\partial_X,\partial_Y)\ \Dalph = 
    \sum_{\bdelta} c_{\bdelta}\,\Dgenalph{\balpha-\bdelta}
  \end{equation}
  for some constants $c_{\bdelta}\in\bbN$.  Here, the sum is over all
  distinct permutations $\bdelta$ of $\bbeta$.  We use the convention
  that $c_{\bdelta}=0$ if $\alpha_i-\delta_{i}\notin \A$ for some $1
  \le i \le n$.
\end{lemma}
\begin{proof}
  First note that $m_{\bbeta}(\partial_X,\partial_Y)$ can be written as 
  \begin{equation*}
    m_{\bbeta}(\partial_X,\partial_Y) = 
    \sum_{\bdelta} \prod_{i=1}^n \partial_{x_i}^{\delta_{i,1}} 
    \partial_{y_i}^{\delta_{i,2}} =  
    K \sum_{\nu\in S_n} \prod_{i=1}^n \partial_{x_i}^{\beta_{\nu(i),1}} 
    \partial_{y_i}^{\beta_{\nu(i),2}}
  \end{equation*}
  for some constant $K\in\mathbb{Q}$ dependent on the extent to which
  factors in $\bbeta$ are repeated.  Then,
  \begin{equation}
  \begin{aligned}
    m_{\bbeta}(\partial_X,\partial_Y)\ \Dalph &=    
    \sum_{\nu\in S_n} K \prod_{i=1}^n \partial_{z_i}^{\beta_{\nu(i)}} 
    \sum_{\sigma\in S_n} \sgn(\sigma) z_1^{\alpha_{\sigma^{-1}(1)}}\cdots 
    z_n^{\alpha_{\sigma^{-1}(n)}}\\
    &= \sum_{\nu\in S_n} c_{\nu}
    \sum_{\sigma\in S_n} \sgn(\sigma) 
    \prod_{i=1}^n z_i^{\alpha_{\sigma^{-1}(i)}-\beta_{\nu(i)}}\\
    &= \sum_{\phi\in S_n} c_{\phi}
    \sum_{\sigma\in S_n} \sgn(\sigma) 
    \prod_{i=1}^n z_i^{\alpha_{\sigma^{-1}(i)}-\beta_{\sigma^{-1}\phi(i)}}\\
    &= \sum_{\bdelta} c_{\bdelta}\, \Dgenalph{\balpha-\bdelta}
  \end{aligned}
  \end{equation}
  for some constants $c_{\bdelta}$.  In the above, we consider the
  coefficient $c_\nu$ to be zero if the exponent of any of the $z_i$'s
  is not in $\A$.  In addition, we have $\nu = \sigma^{-1}\phi$ and
  let $\bdelta$ run over all distinct permutations of $\bbeta$.
\end{proof}

\section{The Ideal $\GS$ and the ring $\RG$.}
Let $\GS$ be the ideal generated by the monomials in the collections 
\begin{align}
  &\left\{x_1y_1, \ldots,x_ny_n \right\},\label{eq:GSone}\\ 
  &\left\{\prod_{i\in D} x_i\right\}_{
    \substack{D \subset \{1, 2, \ldots, n\}\\|D|=m_1+p_1+1}}, \text{ and }
  \left\{\prod_{h\in E} y_h\right\}_{
    \substack{E \subset \{1, 2, \ldots, n\}\\|E|=m_2+p_2+1}}.\label{eq:GStwo}
\end{align}

Note that each of the monomials in equations \eqref{eq:GSone} and \eqref{eq:GStwo}
is in the ideal $\IS$.  Thus $\GS \subset \IS$.

For a sequence $\balpha = (\alpha_1,\ldots,\alpha_n)\in(\A)^n$, we
define 
\begin{equation}
  \begin{aligned}
    |X(\balpha)| &= \alpha_{1,1} + \alpha_{2,1} + \cdots + \alpha_{n,1} \text{ and }\\
    |Y(\balpha)| &= \alpha_{1,2} + \alpha_{2,2} + \cdots + \alpha_{n,2}.    
  \end{aligned}
\end{equation}
For an arbitrary bipermanent $b=[T,V]_{\per}$, we write $|X(b)|$ as
shorthand for $|X(\rowseq(V))|$; similarly for $|Y(b)|$.  In addition
we will write
\begin{equation}
  (q)_j = (1-q)(1-q^2)\cdots(1-q^j)\text{ and }
  (t)_j = (1-t)(1-t^2)\cdots(1-t^j)
\end{equation}
for the \emph{rising factorial products}.  The generating function for
the sum
\begin{equation}
  \sum_{\balpha} t^{|X(\balpha)|}\ q^{|Y(\balpha)|},
\end{equation}
subject to the constraints that $\balpha \in \A^n$,
$\alpha_{i} \leq_{\A} \alpha_{i+1}$ for $1 \le i \le n-1$ and
$\alpha_{m_2+p_2+1}=(0,0)$, is given by
\begin{equation}\label{eq:gfalpha}
  \sum_{\balpha} t^{|X(\balpha)|}\ q^{|Y(\balpha)|} = 
  \frac{1}{(t)_{m_1+p_1}}\frac{1}{(q)_{m_2+p_2}}.
\end{equation}
Define 
\begin{equation}\label{E:Bgamma} 
  \BPgam = \left\{ [T,V]_{\per}:\ T\in\stdtabn,\ 
    V\in\cochgalln{\gamma} \right\}. 
\end{equation}
\begin{theorem}\label{T:Nine} 
  The Hilbert series of $\RG$ is given by
\begin{equation}
  {\mathcal H}(\RG)= \frac{1}{(t)_{m_1+p_1}}\frac{1}{(q)_{m_2+p_2}}
  \sum_{b\in\BPgam} t^{|X(b)|}\ q^{|Y(b)|}.
\end{equation}
\end{theorem}

\begin{proof}
  Suppose $(T,V)\in\bitabsetpn$ and $v_{m_2+p_2+1} \neq (0,0)$.  Since
  $V\in\colstr{\A'}$, either each of the monomials in $[T,V]_{\per}$
  has at least $m_1+p_1+1$ distinct $x_i$'s as factors, or each of the
  monomials has at least $m_2+p_2+1$ distinct $y_i$'s as factors.  In
  either case, $[T,V]_{\per}\in\GS$ as can be seen by examining the
  sets in \eqref{eq:GStwo}.  It follows that in looking for a basis of
  $\RG$, we can restrict our attention to those bipermanents for which
  $v_{m_2+p_2+1} = (0,0)$.

  Additionally, if $v_i\not\in\A$ for some $1\leq i\leq n$, then each
  monomial of $[T,V]_{\per}$ is a multiple of $x_jy_j$ for some $j$.
  (Note that $j$ need not equal $i$ as $\std(T)$ need not equal
  $\std(V)$.)  Examination of \eqref{eq:GSone} then shows that
  $[T,V]_{\per}\in\GS$ in this case as well.

  On the other hand, it is easily seen that if $v_{m_2+p_2+1} = (0,0)$
  and $V\in\colstr{\A}$, then $[T,V]_{\per}\not\in \GS$.  It follows
  from Lemma~\ref{T:lemfour} that the set
  \begin{equation}\label{eq:bppbasis}
    \left\{ \left[T,V\right]_{\per}: 
      (T,V)\in\bitabsetn,\ v_{m_2+p_2+1} = (0,0)\right\}
  \end{equation}
  is a basis for $\RG$. 

  The theorem then follows by combining the decomposition of
  Lemma~\ref{lem:decomposition} with the basis of \eqref{eq:bppbasis}
  and the generating function of \eqref{eq:gfalpha}.
\end{proof}

Although the above proof constructs a basis of $\RG$, it will be more
useful in the next section to have the basis offered by
Theorem~\ref{T:Thefour}.

\begin{theorem}\label{T:Thefour}
  The collection
  \begin{align}\label{E:thirty}
    \EEBbasis = \left\{\left(e_1^{\epsilon_1}(X_n)\cdots
    e_{m_1+p_1}^{\epsilon_{m_1+p_1}}(X_n)\ e_1^{\beta_1}(Y_n)
    \cdots e_{m_2+p_2}^{\beta_{m_2+p_2}}(Y_n)\right)b : b\in\BPgam\right\},
  \end{align}
  where the $\epsilon_i$ and $\beta_i$ are allowed to run over all
  nonnegative integers, is a basis for $\RG$.
\end{theorem}

\begin{proof}
  For this proof, linearly extend the notation $|X(p)|$ to apply to
  elements $p\in \EEBbasis$.  

  By the proof of Theorem~\ref{T:Nine}, it suffices to consider a
  bideterminant $[T,V]_{\per}$ where $(T,V)\in\bitabsetn$ and
  $v_{m_2+p_2+1} = (0,0)$.  Set $C=\cochg{\gamma}(\std(V))$.  Note
  that $c_{m_2+p_2+1} = (0,0)$ by construction.

  If $c_i = v_i$ for all $1\leq i\leq n$, then
  $V\in\cochgalln{\gamma}$; hence, by definition,
  $[T,V]_{\per}\in\BPgam$.  Assume not.  We consider two
  possibilities.

  Suppose there exists such an index greater than $m_2+p_2+1$.  Choose
  $i$ to be the smallest such index.  Let $U$ denote the tableau for
  which $u_j = v_j$ for $1\leq j < i$ and $u_j = v_j - (1,0)$ for
  $i\leq j\leq n$.  It follows from the proof of
  Lemma~\ref{lem:decomposition} that $U\in\colstrn{\A}$.  Note that
  that $(T,U) >_{\per} (T,V)$.

  Utilizing the identity $e_{n-i+1} = m_{1^{n-i+1}}$, it follows from
  Lemma~\ref{T:Cortwo} that
  \begin{equation}\label{E:expanding}
    e_{n-i+1}(X_n)\ [T,U]_{\per} \equiv c_{T,V}\,[T,V]_{\per}  + 
    \sumsb{(S,W) >_{\per} (T,V)\\S\in\stdtabn,\ W\in\colstrn{\A}} 
    c_{S,W}\,[S,W]_{\per} \mod(\GSnxy).
  \end{equation}
  (Here we have used the fact that the monomials in $e_{n-i+1}(X_n)\
  [T,U]_{\per}$ arising in the third term in \eqref{eq:mbeta} all have
  a factor $x_j y_j$ for some $j$.  But, as we see from
  \eqref{eq:GSone}, these monomials are in $\GSnxy$ by construction.)
  The only remaining possibility is that there exists such an index
  $i$ less than $m_2+p_2+1$.  Choose $i$ to be the largest such $i$.
  Arguing as above, we obtain an equivalent expansion for $e_i(Y_n)$.

  Iteration of the above argument on the $[T,U]_{\per}$ and
  $[S,W]_{\per}$ implies that the collection $\EEBbasis$
  spans $\RG$.  (We use the facts that at most $m_1+p_1$ of the $x_i$,
  and $m_2+p_2$ of the $y_i$, can appear with positive degree in any
  of the monomials not in $\GSnxy$.)

  For any $b\in \BPgam$, let $\pow(\bepsilon,X_n,b)$ denote the
  $X_n$-degree of $\prod_{i=1}^{m_1+p_1} e_i^{\epsilon_i}(X_n) b$ and
  $\pow(\bbeta,Y_n,b)$ denote the $Y_n$-degree of $\prod_{i=1}^{m_2+p_2}
  e_i^{\beta_i}(Y_n) b$.  It follows then that 
  \begin{equation}
    \begin{aligned}
      \sum_{p \in \EEBbasis} t^{|X(p)|}\ q^{|Y(p)|} 
      &= \sum_{\bepsilon\in\bbN^{m_1+p_1}} \sum_{\bbeta\in\bbN^{m_2+p_2}}
      \sum_{b \in \BPgam} t^{\pow(\bepsilon,X_n,b)} q^{\pow(\bbeta,Y_n,b)}\\ 
      &= \frac{1}{(t)_{m_1+p_1}}\frac{1}{(q)_{m_2+p_2}}
      \sum_{b\in\BPgam} t^{|X(b)|} q^{|Y(b)|}.
    \end{aligned}
  \end{equation}
  The fact that the collection $\EEBbasis$ spans $\RG$ and yields
  the desired Hilbert series (see Theorem \ref{T:Nine}) implies that
  $\EEBbasis$ must be a basis for $\RG$.
\end{proof}

We have the following corollary:
\begin{corollary}\label{T:Corfive}
  The collection
  \begin{equation}\label{eq:EEB}
    \left\{e_1(X_n), e_2(X_n), \ldots, e_{m_1+p_1}(X_n),
    e_1(Y_n), e_2(Y_n),\ldots, e_{m_2+p_2}(Y_n)\right\}
  \end{equation}
  is algebraically independent in the ring $\RG.$
\end{corollary}
\begin{proof}
  Any nontrivial algebraic dependence amongst the elements of
  \eqref{eq:EEB} would yield a linear dependence amongst the elements
  of $\EEBbasis$, conflicting with its role in the basis for $\RG$.
\end{proof}

\section{The ideal $\HS$ and the ring $\RH$.}
\label{sec:HR}

Let $\HS$ be the ideal in $\CXY$ generated by the collections of
monomials in equations \eqref{eq:GSone} and \eqref{eq:GStwo} as well
as by the elementary symmetric polynomials in the collection
\begin{equation}\label{eq:hdef}
  \left\{e_{p_1+2}(X_n),\ldots, e_{p_1+m_1}(X_n), 
  e_{p_2+2}(Y_n),\ldots, e_{p_2+m_2}(Y_n)\right\}.
\end{equation}
\myfig{.4}{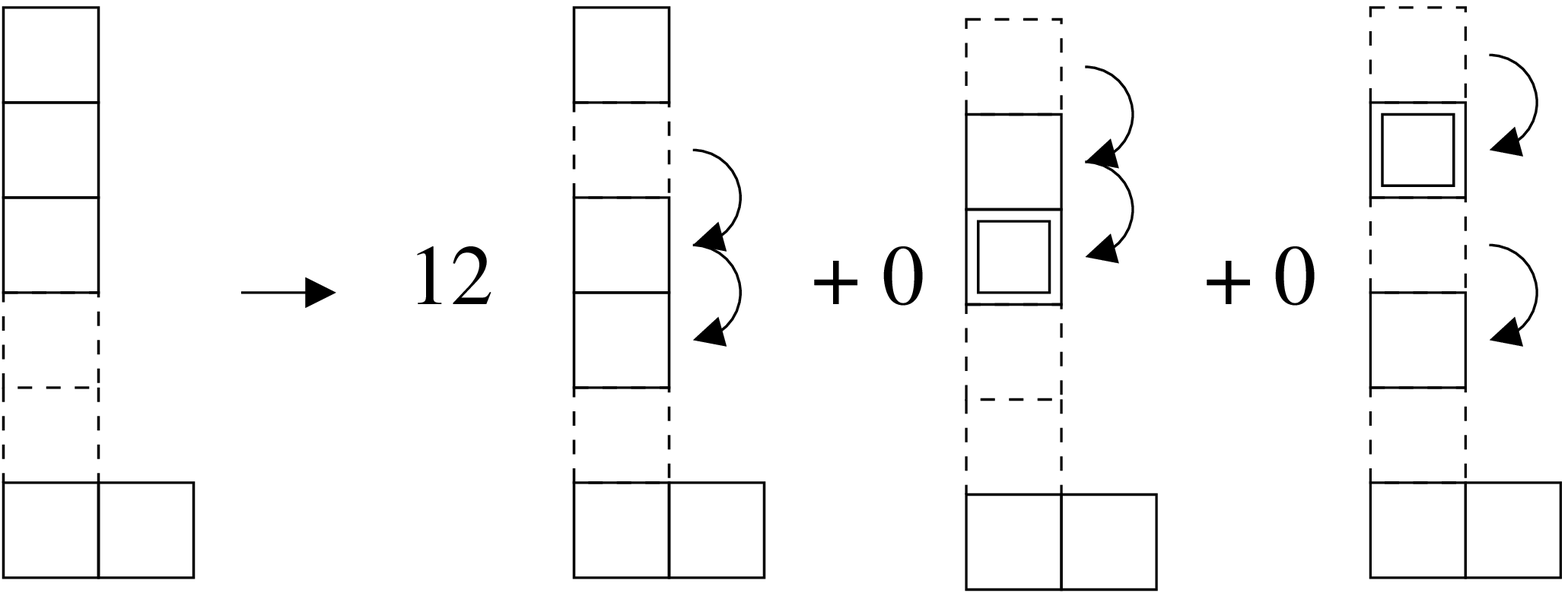}{Illustration of the fact that
  $e_2(\partial_Y)\Dgenalph{(\olf,\ole,\old,\ula,\ulb)}$ equals
  $12\Dgenalph{(\olf,\old,\olc,\ula,\ulb)}$.}

Consider the action of an elementary symmetric differential operator
$e_j(\partial_{Y_n})$, such as is illustrated in
Figure~\ref{fig:eaction} for $j=2$.  This operator moves each of $j$
distinct cells down by one place.  Any configuration in which two
cells end up in the same position or in which a cell moves to a
position not indexed by an element of $\A$ contributes zero.  (The
action of an $e_l(\partial_{X_n})$ is similar.)  It follows that for
any contributing monomial, the cells contiguous with $(0,0)$ are not
moved.  But any $e_{p_1+j}(\partial_X)$ for $j > 1$ or
$e_{p_2+\ell}(\partial_Y)$ for $\ell > 1$ \emph{must} move one of these
fixed cells.  So $\HS \subset \IS$.  Note that by construction, $\HS
\supset \GS$.  Theorem \ref{T:Thefour} and Corollary \ref{T:Corfive}
imply the following two corollaries:
\begin{corollary} 
  The set
  \begin{equation}
    \left\{ \left(\prod_{i=1}^{p_1+1} e_i^{\epsilon_i}(X_n)\cdot 
            \prod_{j=1}^{p_2+1} e_j^{\epsilon_j}(Y_n)\right) b : 
            b\in\BPgam\right\}_{\substack{\bepsilon\in\bbN^{p_1+1}\\
                                          \bbeta\in\bbN^{p_2+1}}},
\end{equation}
is a basis for $\RH$.
\end{corollary}

\begin{corollary}
  The Hilbert series of $\RH$ is given by
  \begin{equation}
    {\mathcal H}(\RH) = \frac{1}{(t)_{p_1+1}}\frac{1}{(q)_{p_2+1}}
    \sum_{b\in\BPgam} t^{|X(b)|} q^{|Y(b)|}.
  \end{equation}
\end{corollary}

\section{The ideals $\JS$ and $\KS$}

In the previous section, we considered symmetric functions whose
corresponding differential operators moved collections of boxes, but
for which each cell was only moved one place.  Monomial symmetric
functions yield operators that move cells farther.  In this section we
consider which ones will also annihilate $\Dalph$.  In fact, due to
fortuitous cancellations, we will focus on the differential operators
corresponding to the complete symmetric functions.

\myfig{.4}{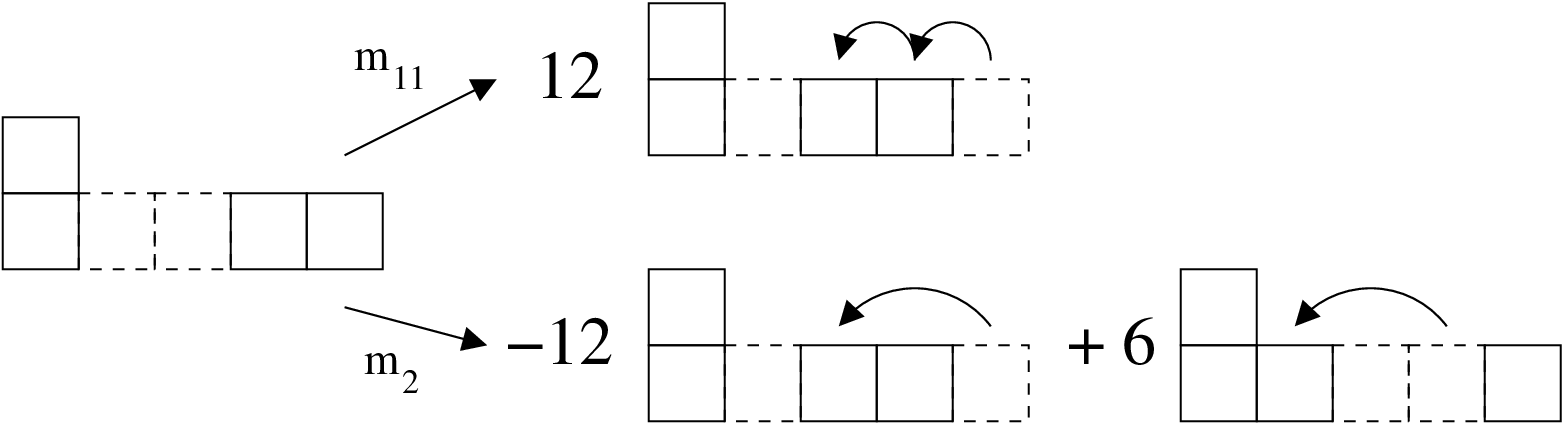}{Illustration of the action of $h_2(\partial_{X_4})
  = m_2(\partial_{X_4}) + m_{11}(\partial_{X_4})$ on $\Dalph$ for
  $\balpha = ((0,1),(0,0),(3,0),(4,0))$.}

Figure~\ref{fig:maction} illustrates the action of $h_2(\partial_X) =
m_2(\partial_X) + m_{11}(\partial_X)$ on a given $\Dalph$.  Notice
that the hollow lattice diagram $\tilde{\gamma} = ((1,2),(2,0),(1,0))$
is obtained in two different ways: once through the action of
$m_{11}(\partial_X)$ and once through $m_2(\partial_X)$.  However, for
$m_2(\partial_X)$ the cell that moves jumps over another cell.  This
leads to the introduction of a sign.  Hence, the two
$\Delta_{\tilde{\gamma}}$ that appear cancel.  In fact, as
Lemma~\ref{lem:inv} shows, under the action of an $h_j$ on some
$\Dalph$, the only term that survives is that which moves the cell
$(m_1+k_1-1,0)$ (or $(0,m_2+k_2-1)$, as appropriate) $j$ spaces.

We now need to consider lattice diagrams that are subsets of hook
shapes, but are not hollow.  In particular, we wish to consider hollow
diagrams modified by sliding certain cells closer to the origin.  The
amount of sliding will be described by two sequences $a_0,\ldots,a_i$
and $b_0,\ldots,b_j$ of nonincreasing, nonnegative integers.  For
brevity in what follows, write
\begin{equation*}
  \boldsymbol{c} = (\ol{m_2-1},\ldots,\ol{1},\ul{0},\ul{1},\ldots,\ul{m_1-1}).
\end{equation*}
Then, for $\hollatdiag = \latdiag{\balpha}$ with
\begin{equation*}
  \balpha = (\ol{m_2+k_2+p_2-1},\ldots,\ol{m_2+k_2-1},\boldsymbol{c},
  \ul{m_1+k_1-1},\ldots,\ul{m_1+k_1+p_1-1}),
\end{equation*}
we write $\gamma[a_0,\ldots,a_i; b_0,\ldots,b_j]$ for the lattice
diagram corresponding to the collection
\begin{gather*}
  (\ol{m_2+k_2-1+p_2},\ldots,\ol{m_2+k_2-1+i-a_i},
  \ldots,\ol{m_2+k_2-1+0-a_0},\boldsymbol{c},\\
  \ul{m_1+k_1-1+0-b_0},\ldots,\ul{m_1+k_1-1+j-b_j},\ldots,\ul{m_1+k_1-1+p_1}).
\end{gather*}
\myfig{.4}{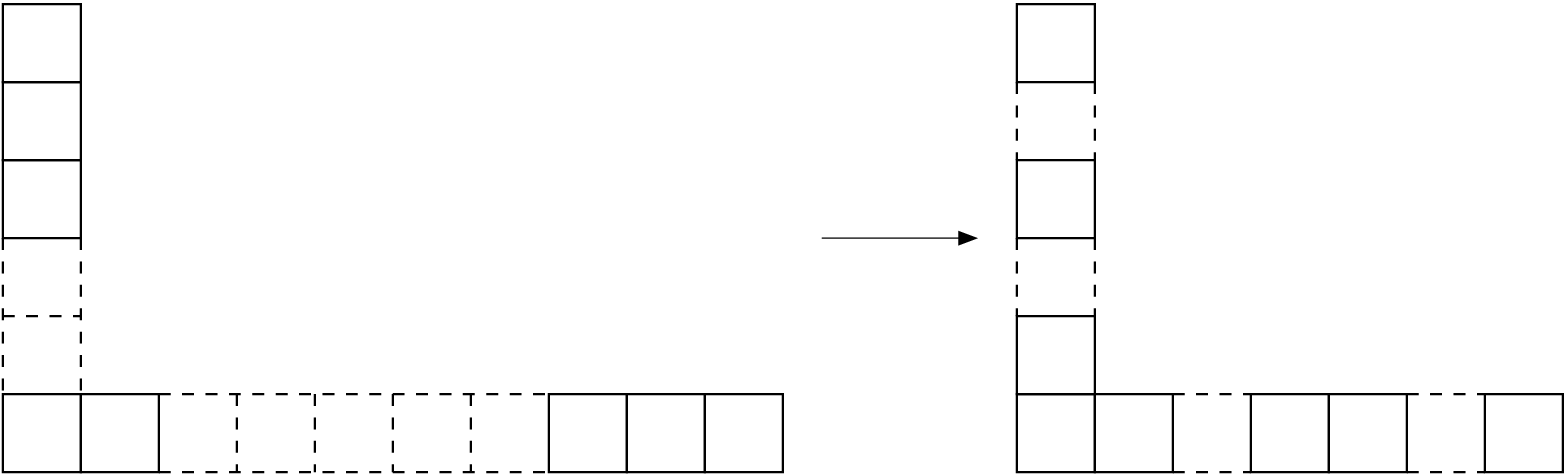}{Sample square bracket notation for $\gamma$.}
Note that due to the nonincreasing restriction on the sequences, 
\begin{equation*}
  \sgn(\DSnoxy) = \sgn(\Delta_{\gamma[a_0,\ldots,a_i; b_0,\ldots, b_j]}).
\end{equation*}

An example is illustrated in Figure~\ref{fig:brack}.  In the figure,
the left diagram is the hollow lattice diagram for $\gamma =
((2,1),(6,3),(2,2))$, while the right diagram illustrates $\gamma[2,1;
4,4,3]$.

\begin{lemma}\label{lem:inv}
  Let $0\leq j<k_1$ and $0\leq \ell < k_2$.  Then
  \begin{equation}\label{E:lowerthem}
    h_j(\partial_X)\ \DS = c_j\ \Delta_{\gamma[0;j]}
    \quad\text{ and }\quad
    h_\ell(\partial_Y)\ \DS = c_\ell\ \Delta_{\gamma[\ell;0]}.
  \end{equation}
\end{lemma}
  
\begin{proof}
  We only prove the $h_j(\partial_X)$ version as the proof
  of the $h_\ell(\partial_Y)$ version is effectively identical.

  Recall that $h_j(X_n) = \sum_{|\la|=j} m_{\la}(X_n)$.  Let us
  consider $\DS$ for $\hollatdiag = \latdiag{\balpha}$ where
  \begin{multline}\label{eq:alphaorder}
    \balpha = \left(\ul{0},\ul{1},\ldots,\ul{m_1-1},
      \ul{m_1+k_1-1},\ul{m_1+k_1+p_1-1},\right.\\
    \left. \ol{1},\ldots,\ol{m_2-1},\ol{m_2+k_2-1},\ldots,\ol{m_2+k_2+p_2-1}\right).
  \end{multline}
  (Note the unusual order in which we have listed the cells.)  Let
  $\lambda = (\lambda_1,\ldots,\lambda_n) \vdash j$.  View $\lambda$
  as an element $\blambda$ of $\A^n$ under the map $\lambda_i \mapsto
  (\lambda_i,0)$.  It follows from Lemma~\ref{T:expands} that
  \begin{equation}\label{E:uglyformthree}
    m_{\blambda}(\partial_X)\ \DS  = \sum_{\bnu} c_{\bnu}\ \Dgenalph{\balpha-\bnu},
  \end{equation}
  where the sum is over all distinct permutations $\bnu$ of
  $\blambda$.  As previously, we use the convention $c_{\bnu} = 0$ if
  $\balpha-\bnu\not\in \A^n$.  

  Consider a particular permutation $\bnu$ of $\blambda$.  If
  $\alpha_i-\nu_{i} = \alpha_j-\nu_{j}$ for some $i \ne j$ then
  $\Dgenalph{\balpha-\bnu}=0$ as $\Dgenalph{\balpha-\bnu}$ is a
  determinant.  Recalling the ordering of $\balpha$ given in
  \eqref{eq:alphaorder}, note that if $\nu_{i,1}\ne 0$ for some $1 \le
  i \le m_1$, then we must have $c_{\bnu}=0$.  Therefore, without loss
  of generality, we may assume that $\nu_{i,1} = 0$ for $1 \le i \le
  m_1$ and that $\alpha_i-\nu_i \neq \alpha_j-\nu_j$ for $i\neq j$.

  \myfig{.4}{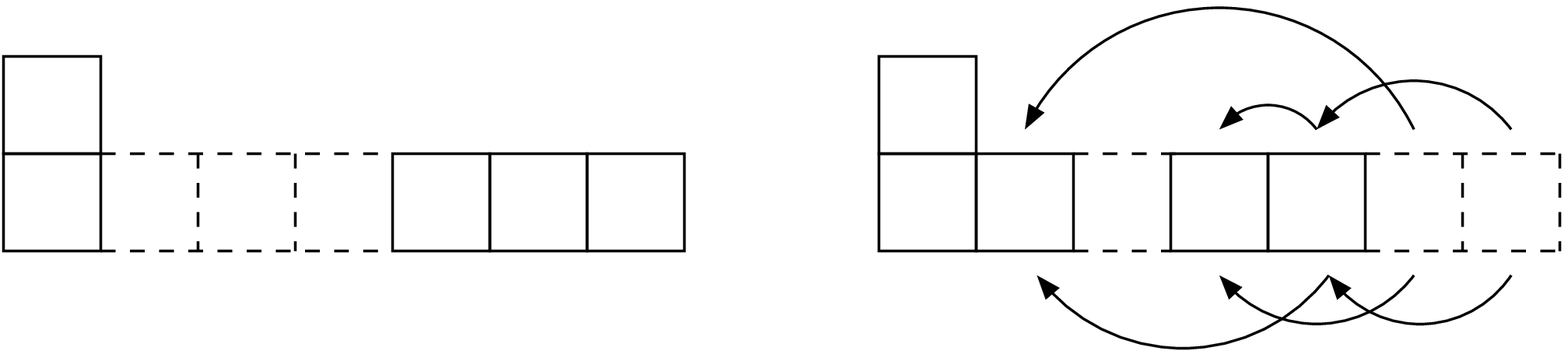}{Illustration of the involution described in the
    proof of Lemma~\ref{lem:inv}.}

  Let $\bmu$ be given by $\mu_{m_1+1} = (j,0)$ and $\mu_i = (0,0)$ for
  $i\neq m_1+1$.  We now proceed to define an involution on the $\bnu$
  whose only fixed point is $\bmu$.  So assume $\bnu \neq \bmu$.

  Define $p < r$ to index the smallest two terms (with respect to
  $<_{\A}$) of
  \begin{equation}
    \alpha_{m_1+1}-\nu_{m_1+1},\alpha_{m_1+2}-\nu_{m_1+2},
    \ldots,\alpha_{m_1+p_1+1}-\nu_{m_1+p_1+1}.
  \end{equation}
  The indices $p$ and $r$ index the two boxes of $\Dalph$ not
  contiguous with the cell $(0,0)$ that have moved farthest to the
  left upon subtraction of $\bnu$.  Set
  \begin{equation}
    q = \left(\alpha_{r,1}-\nu_{r,1}\right)-\left(\alpha_{p,1}-\nu_{p,1}\right).
  \end{equation}
  We now define
  \begin{equation}
    \bnu' = g(\bnu) = 
      \left[ \nu_1, \nu_2, \ldots, \ul{\nu_{p,1}-q}, 
        \ldots, \ul{\nu_{r,1}+q},\ldots,\nu_n\right].
  \end{equation}
  The left-hand picture in Figure~\ref{fig:inv} illustrates $\Dalph$
  for $\balpha = (\ul{0},\ul{4},\ul{5},\ul{6},\ol{1})$.  The right-hand
  picture illustrates $\Dgenalph{\balpha-\bnu}$ for $\bnu =
  (\ula,\ul{3},\ul{2},\ul{2},\ula)$ via the bottom triple of arrows
  ($p = 2$, $r = 3$, $\ul{1} = \alpha_2-\nu_2 < \alpha_3-\nu_3 =
  \ul{3}$, and $q = 2$) and $\Dgenalph{\balpha-\bnu'}$ for $\bnu' =
  (\ula,\ul{1},\ul{4},\ul{2},\ula)$ via the top triple of arrows.
  When viewed as sets, $\balpha-\bnu$ equals $\balpha-\bnu'$; they
  differ only in order.

  We have that
  \begin{align}
    \Dgenalph{\balpha-\bnu} = -\Dgenalph{\balpha-\bnu'},
  \end{align}
  since $\balpha-\bnu$ and $\balpha-\bnu'$ differ by a transposition.
  Now $g(g(\bnu))=\bnu$.  As desired, this function $g$ yields a
  sign-reversing involution between all the terms
  $\Dgenalph{\balpha-\bnu}$ in equation \eqref{E:uglyformthree} except
  for the unique $\Dgenalph{\balpha-\bmu}$. It is not difficult to see
  that the corresponding coefficients $c_{\bnu'}$ in the expansion of
  $m_\lambda(\partial_X) \DS$ in terms of the
  $\Dgenalph{\balpha-\bnu}$ satisfy $c_{\bnu} = c_{\bnu'}$.
  Thus, all the terms in equation \eqref{E:uglyformthree} cancel out except
  $\Dgenalph{\balpha-\bmu}$.  This gives the lemma.
\end{proof}

Recall that $\DS=0$ if two of the entries in $\hollatdiag$ are identical.
Thus, if $j \ge k_1$ or $\ell \geq k_2$, then $\Delta_{\gamma[0;j]} =
0$ or $\Delta_{\gamma[\ell;0]} = 0$, respectively.  It follows that
\begin{corollary}\label{T:corfifteen}
  $h_{k_1+i}(X_n)\in\IS$ for $i \ge 0$ and $h_{k_2+h}(Y_n)\in\IS$ for $h
  \ge 0$. 
\end{corollary}

Corollary~\ref{T:corfifteen} provides information about certain
elements that must be in the ideal $\IS$. As such, we will use it to
define a sub-ideal $\KS$ of $\IS$. Specifically, set $\KS$ to be the
ideal in $\CXY$ generated by the generators of $\HS$ along with
\begin{equation}\label{E:smallhsxys}
  \left(h_{k_1}(X_n), \ldots, h_{k_1+p_1+1}(X_n),\ldots,
    h_{k_2}(Y_n), \ldots, h_{k_2+p_2+1}(Y_n),\ldots\right).
\end{equation}

As it turns out, $\KS$ is a finitely generated ideal in $\RH$.
To this end, define $\JS$ to be the sub-ideal of $\KS$ generated by
the generators of $\HS$ along with 
\begin{equation}\label{E:smallerhsxys}
  \left\{h_{k_1}(X_n), \ldots, h_{k_1+p_1}(X_n), 
    h_{k_2}(Y_n), \ldots, h_{k_2+p_2}(Y_n)\right\}.
\end{equation}

It follows from Corollary~\ref{T:corfifteen} that $\JSnxy \subset
\KSnxy \subset \ISnxy$.  
\begin{lemma}
  $\KS \equiv \JS \pmod \HS$.
\end{lemma}
\begin{proof}
  A standard result (cf. \cite[pg. 21]{Macdonald}) is that
  \begin{equation}
    \sum_{r=0}^n (-1)^r\ e_r(X_n)\ h_{n-r}(X_n)  =  0.  
  \end{equation}
  We prove that $h_{k_1+p_1+a}(X_n)\equiv 0 \pmod \JS$ for $a\geq 1$
  by induction on $a$.  For $a\ge 1$, we have
  \begin{multline}
    h_{k_1+p_1+a}(X_n)
    = \sum_{r=1}^{p_1+1}(-1)^{r+1}\
    e_r(X_n)\ h_{k_1+p_1+ a-r}(X_n)\\ 
    + \sum_{r=p_1+2}^{k_1+p_1+ a}(-1)^{r+1}\ e_r(X_n)\
    h_{k_1+p_1+ a-r}(X_n).\label{eq:recurIII}
  \end{multline}
  Recall from \eqref{eq:GStwo} and \eqref{eq:hdef} that
  $e_{p_1+j}(X_n)\in\HS$ for $j\ge 2$.  So the second sum of
  \eqref{eq:recurIII} is in $\HS\subset \JS$.  On the other hand, by
  the definition of $\JS$, $h_{k_1+i}\in\JS$ for $0\leq i\leq p_1$.
  So the first sum of \eqref{eq:recurIII} is in $\JS$ as well.
  This proves the claim for $a = 1$.  The claim for $a > 1$ thereby
  follows by the obvious induction hypothesis.  Similar arguments can
  be made about $h_{k_2+p_2+ a}(Y_n)$, for $ a \ge 1$.
\end{proof}

Our next goal is to show that the collection that generates $\JS$ in
equation \eqref{E:smallerhsxys} is itself algebraically independent in
$\RH$.

\begin{theorem}\label{T:funtheorem}
  The collection
  \begin{equation}
    \left\{h_{k_1}(X_n), h_{k_1+1}(X_n),
    \ldots, h_{k_1+p_1}(X_n), h_{k_2}(Y_n), h_{k_2+1}(Y_n), \ldots,
    h_{k_2+p_2}(Y_n)\right\}
  \end{equation}
  is algebraically independent in the ring $\RH$.
\end{theorem}

The reader is advised to follow Example~\ref{ex:he} while reading the proof.
\begin{proof}
  The theorem is equivalent to the statement that there is no
  nontrivial polynomial $P$ over $\bbC$ in $p_1+p_2+2$ variables such that 
  \begin{equation}
    P(h_{k_1}(X_n),\ldots,h_{k_1+p_1}(X_n),h_{k_2}(Y_n),\ldots,h_{k_2+p_2}(Y_n))
    \in \HS.
  \end{equation}
  Since the $x_i$ and $y_i$ are true indeterminates, such a relation
  would have to continue to hold upon the specialization $y_1 = \cdots
  = y_n = 0$.  Hence, it suffices to show that there is no nontrivial
  polynomial $Q$ such that
  \begin{equation}\label{eq:algdep}
    Q(h_{k_1}(X_n),\ldots,h_{k_1+p_1}(X_n)) \in \HS.
  \end{equation}
  So, to prove the theorem, we assume that such a $Q$ does exist and
  obtain a contradiction.  Specifically, we will show that if the
  $h_i$ in question are algebraically dependent in $\HS$, then the
  $e_i$ of \eqref{eq:EEB} are algebraically dependent in $\GS$, in
  direct contradiction with Lemma~\ref{T:Corfive}.

  A relation such as \eqref{eq:algdep} can be rewritten as 
  \begin{equation}\label{eq:hsum}
    \sum_\la c_\la\ h_\la \in \HS,
  \end{equation}
  where each $\la$ is a partition with parts of length chosen from the
  collection $\{k_1, k_1+1, \ldots, k_1+p_1\}$.  We assume that $c_\la
  = 0$ for any $\la$ with $h_\la\in \HS$ and that there exists some
  $c_\la \neq 0$.

  Since the monomial symmetric functions are a basis of the ring of
  symmetric functions, any such $h_\la$ can be expanded in terms of
  the $e_\mu$.  In fact, this can be done explicitly as follows (cf.
  \cite[pg. 107]{Macdonald}).

  Define a \emph{domino} to be a set of horizontally consecutive
  squares in a Ferrers diagram.  For $\la,\mu\vdash n$, a \emph{domino
    tabloid of shape $\la$ and type $\mu$} is a tiling of
  $\diagram(\la)$ with dominoes of length $\mu_1,\mu_2,\ldots$.  We
  consider dominoes of the same length to be indistinguishable.  Let
  $d_{\la,\mu}$ denotes the number of domino tabloids of shape $\la$
  and type $\mu$.  Then we can write
  \begin{equation}\label{eq:domexpn}
    h_\la (X_n) =  \sum_{\mu} (-1)^{|\mu|-\ell(\mu)} d_{\la,\mu}\ e_\mu (X_n).
  \end{equation}
  Recall that $e_\mu$ is in $\GS$ (respectively, $\HS$) whenever $\mu$
  has a part greater than or equal to $m_1+p_1+1$ (respectively,
  $p_1+2$).  Combining \eqref{eq:domexpn} and \eqref{eq:hsum}, we find
  that
  \begin{equation}
    \sumsb{\mu\\\mu_1\leq p_1+1} (-1)^{|\mu|-\ell(\mu)} \left(\sum_\lambda c_\lambda 
    d_{\la,\mu}\right) e_\mu (X_n) \equiv 0 \pmod{\HS},
  \end{equation}
  or, equivalently, that
  \begin{multline}\label{eq:hedep}
    \sumsb{\mu\\\mu_1\leq p_1+1} (-1)^{|\mu|-\ell(\mu)} \left(\sum_\lambda c_\lambda 
    d_{\la,\mu}\right) e_\mu (X_n)\\
    - \sumsb{\mu\\ p_1+2\leq \mu_1\leq p_1+m_1} a_\mu e_\mu (X_n) \equiv 0 \pmod{\GS}.
  \end{multline}
  for some constants $a_\mu$.  If one of the $a_\mu \neq 0$ or one of
  the $\sum_\lambda c_\lambda d_{\la,\mu} \neq 0$ for some $\mu$, then
  we have a contradiction (asserted by Corollary~\ref{T:Corfive}) of the
  algebraic independence in $\RG$ of the set $\{e_i\}_{1\leq i\leq
    p_1+m_1}$.  We will, in fact, show that at least one of the
  $\sum_\lambda c_\lambda d_{\la,\mu}$ is nonzero.

  For a given $\lambda$, there is a maximum (with respect to
  lexicographic order) $\mu$ with $d_{\la,\mu} \neq 0$ and $e_\mu\not
  \in \HS$.  It is given by
  \begin{equation}\label{E:type}
    L(\la) =  \left( (p_1+1)^{\alpha_{p_1+1}},
      (p_1)^{\alpha_{p_1}}, (p_1-1)^{\alpha_{p_1-1}}, \ldots,
      (1)^{\alpha_{1}} \right). 
  \end{equation}
  Here, for $1 \le f \le p_1$, $\alpha_{f}$ equals the number of
  entries in
 \begin{equation}
   \la = (\la_1,\la_2,\ldots,\la_{\ell(\lambda)})
 \end{equation}
 congruent to $f \pmod{p_1+1}$
 and
 \begin{equation}
   \alpha_{p_1+1} =  \sum_{r=1}^{\ell(\lambda)} \left\lfloor 
     \frac{\la_r}{p_1+1}\right\rfloor
 \end{equation}
 ($\lfloor x \rfloor$ denotes the {\it floor of $x$}).  In particular,
 modulo $\HS$, we can rewrite \eqref{eq:domexpn} as 
  \begin{equation}\label{eq:domexpnII}
    h_\la (X_n) \equiv \pm d_{\la,L(\la)} e_{L(\la)}(X_n) + 
    \sum_{\mu \revlexordstr L(\la)} (-1)^{|\mu|-\ell(\mu)} d_{\la,\mu}\ e_\mu (X_n)
  \end{equation}
  with $d_{\la,L(\la)} \neq 0$.

  Furthermore, given $\mu$, there is a unique partition $\la$
  (possibly with $c_\la = 0$) for which $\mu$ is the maximum element
  appearing in \eqref{eq:domexpn} with $d_{\la,\mu} \neq 0$ for some
  $\tilde{\lambda}$.  To see this, note that each $\alpha_f$ (when $0 \le f
  \le p_1$) gives the number of $\tildlam_j$ congruent to $f
  \pmod{p_1+1}$.  Since $k_1 \le \tildlam_j \le k_1 + p_1$ this
  uniquely identifies $\tildlam_j$.  Each such $\tildlam_j$ requires
  $\left\lfloor \frac{\tildlam_j}{p_1+1}\right\rfloor$ dominoes of
  length $p_1+1$ to finish out the row.  Any remaining $p_1+1$ will be
  used to construct the remaining rows of $\tildlam$.
 
  To finish the proof, pick from \eqref{eq:hsum} the partition
  $\tildlam$ occurring with $c_{\tildlam} \neq 0$ for which
  $L(\tildlam)$ is maximal.  Certainly $d_{\tildlam,L(\tildlam)} \neq 0$.
  However, by this choice of $\tildlam$, $d_{\la,L(\tildlam)} = 0$ or
  $c_\la = 0$ for $\la \neq \tildlam$.  Hence, the coefficient of
  $e_{L(\tildlam)}$ in \eqref{eq:hedep} is nonzero as desired.
\end{proof}

\begin{example}\label{ex:he}
  Fix $\gamma = ((4,1),(7,0),(2,0))$.  In this example we consider the
  expansion of $h_\lambda = h_{(9,9,9,8,7,7)}$ in terms of the
  monomial symmetric functions modulo $\HS$.  In particular, we
  construct $L(\la)$ and show that $\la$ is recoverable from it.

  From $p_1 = 2$ it follows that $L((9^3,8,7^2)) = (3^{15},2,1^2)$.
  The dominoes on each of the rows of length $7$ and $8$ can be
  placed in three different ways, so we find
  $d_{(9^3,8,7^2),(3^{15},2,1^2)} = 27$.  Then, following
  \eqref{eq:domexpnII}, $h_{(9^3,8,7^2)}(X_n)$ is congruent to
  \begin{equation}
     -27 e_{(3^{15},2,1^2)}(X_n) + \\
    \sum_{\mu <_{\normlexord} (3^{15},2,1^2)} (-1)^{|\mu|-\ell(\mu)} 
    d_{(9^3,8,7^2),\mu} e_{\mu}(X_n)
  \end{equation}
  modulo $\HS$.  On the other hand, given $L(\la) = (3^{15},2,1^2)$,
  we can recover $\la$ as follows.  We know that the parts of $\la$
  must all be between $k_1$ and $k_1+p_1$; $7$ and $9$ in this
  case.  The two $1$'s in $L(\la)$ tell us that there must be two
  parts of $\la$ of length $1$ modulo $p_1 + 1 = 3$ (i.e., two rows of
  length $7$).  Similarly, we compute that there is a unique row of
  length $8$.  These three rows account for six of the fifteen length
  $3$ parts of $L(\la)$.  The remaining nine length $3$ parts must
  together comprise the parts of length $9$.  Hence there are three of
  them and we recover $\la = (9^3,8,7^2)$ as desired.  This completes
  our example.
\end{example}

Theorem \ref{T:funtheorem} implies the following.

\begin{corollary}\label{cor:final}
  The Hilbert Series of $\RJ$ is given by

  \begin{equation}
\begin{aligned}\label{E:HilbJ} 
  \mathcal{H}(\RJ) &= \frac{(t)_{k_1+p_1}}{(t)_{k_1-1}(t)_{p_1+1}}
                      \frac{(q)_{k_2+p_2}}{(t)_{k_2-1}(t)_{p_2+1}}
                      \sum_{b\in\BPgam} t^{|X(b)|} q^{|Y(b)|}\\  
                   &= \left[
    \begin{matrix}
      p_1+k_1\\
      p_1+1\cr
    \end{matrix}
  \right]_t \left[
    \begin{matrix}
      p_2+k_2\\ p_2+1\cr
    \end{matrix}
  \right]_q \sum_{b\in\BPgam} t^{|X(b)|}\ q^{|Y(b)|}.
\end{aligned}
  \end{equation}
\end{corollary}

\section{The ideals $\JS = \IS$.}

We need to identify the generators of $\IS$ (recall equation
\eqref{E:defIS}). Specifically, we want to prove that it is finitely
generated by a collection of complete symmetric functions in the ring
$\RH$. The goal is to show that the ideals $\JSnxy$ and $\ISnxy$ are
equal in $\RH$.  Since $\JSnxy\subseteq\ISnxy$, we will do this by
constructing a linearly independent set in $\RI$ that gives the
Hilbert series in equation \eqref{E:HilbJ}.  Set
\begin{equation}
  Q_{k,p} = \left\{(q_1, \ldots, q_p)\in\bbN^p:
    \text{ the } q_i \text{ are nonincreasing and } q_1\leq k\right\}.
\end{equation}

Observe that $|Q_{k,p}| = \binom{p + k}{p}$.  We can now define the
collection of polynomials $\hhBPgam$ that will turn out to be the
required basis for $\RI.$ (Recall that $\BPgam$ is defined in equation
\eqref{E:Bgamma}.)
\begin{equation}\label{D:basisdef}
  \hhBPgam = \left\{ h_q(X_n) h_{q'}(Y_n) \ b :b\in\BPgam,
  q \in Q_{k_1-1,p_1+1}, q'\in Q_{k_2-1,p_2+1} \right\}.
\end{equation}

It follows from the remark above that the number of elements in
$\hhBPgam$ is
\begin{equation}
  n! \binom{p_1+k_1}{p_1+1}\binom{p_2+k_2}{p_2+1}.  
\end{equation}

Furthermore, note that
\begin{equation}\label{E:rightsum} 
  \sum_{b \in \hhBPgam} t^{|X(b)|}\ q^{|Y(b)|}
   =  \left[
    \begin{matrix}
      p_1+k_1\\ p_1+1\cr
    \end{matrix}
  \right]_t \left[
    \begin{matrix}
      p_2+k_2\\ p_2+1\cr
    \end{matrix}
  \right]_q \sum_{b\in\BPgam} t^{|X(b)|}\ q^{|Y(b)|},
\end{equation}
which equals the summation found in equation \eqref{E:HilbJ}.

In Theorem~\ref{thm:hdelta} we will consider the action of
differential operators in the complete symmetric functions on the
determinants $\Delta_\gamma$.  This theorem generalizes Lemma~\ref{lem:inv}.
The below example gives a sample computation in this spirit.
\begin{example}
  By Lemma~\ref{lem:inv}, we have
  \begin{equation}\label{eq:delexpan}
  \begin{aligned}
    h_{3,2}(\partial_X)\
    \Delta_{L[(\ul{0},\ul{1},\ul{2},\ul{3},\ul{9},\ul{10},\ul{11})]}
    &= h_2(\partial_X)\ h_3(\partial_X)\ 
    \Delta_{L[(\ul{0},\ul{1},\ul{2},\ul{3},\ul{9},\ul{10},\ul{11})]} \\ 
    &= \left((m_2 + m_{1,1})(\partial_X)\right)
    (9\cdot 8\cdot 7)\Delta_{L[(\ul{0},\ul{1},\ul{2},\ul{3},\ul{6},\ul{10},\ul{11})]}
  \end{aligned}
  \end{equation}
  equal to
  \begin{equation}\label{eq:delexpanII}
  \begin{gathered}
     \frac{11!}{6!}\
     \Delta_{L[(\ul{0},\ul{1},\ul{2},\ul{3},\ul{6},\ul{10},\ul{9})]}
     + \frac{11!}{6!}\
     \Delta_{L[(\ul{0},\ul{1},\ul{2},\ul{3},\ul{6},\ul{9},\ul{10})]}
     + 9\frac{10!}{6!}\
     \Delta_{L[(\ul{0},\ul{1},\ul{2},\ul{3},\ul{6},\ul{8},\ul{11})]}\\
     + \frac{10!}{5!}\
     \Delta_{L[(\ul{0},\ul{1},\ul{2},\ul{3},\ul{5},\ul{9},\ul{11})]}
     + \frac{9!}{4!}\
     \Delta_{L[(\ul{0},\ul{1},\ul{2},\ul{3},\ul{4},\ul{10},\ul{11})]}.
  \end{gathered}
  \end{equation}
  (We have omitted those $\Dalph$ that have repeated entries and are
  thus identically equal to zero.  We have also included exact
  coefficients even those these are not given by Lemma~\ref{lem:inv}.)
  Note that the first and second terms cancel as they differ by the
  transposition of adjacent elements.  This completes our example.
\end{example}
Suppose the sequences $\bbeta$ and $\bdelta$ correspond to $\gamma' =
\gamma[a_0,\ldots, a_i; b_0,\ldots b_j]$ and $\gamma'' =
\gamma[a_0',\ldots, a_i'; b_0',\ldots b_j']$, respectively, for some
triple $\gamma$ and nonincreasing sequences $(a_\ell)$, $(b_\ell)$,
$(a_\ell')$ and $(b_\ell')$.  We write $\bbeta >_{\content} \bdelta$
(or $\gamma' >_{\content} \gamma''$) whenever
\begin{enumerate}
\item $(a_0,\ldots,a_i) >_{\normlexord} (a_0',\ldots,a_i')$; or
\item if $(a_0,\ldots,a_i) = (a_0',\ldots,a_i')$ then
  $(b_0,\ldots,b_j) >_{\normlexord} (b_0',\ldots,b_j')$.
\end{enumerate}

\begin{theorem}\label{thm:hdelta}
  If 
  \begin{equation}
    \begin{aligned}
      q &= (q_1, q_2, \ldots, q_{p_1+1})\in Q_{k_1-1,p_1+1} \ \text{ and}\\
      q' &= (q_1', q_2', \ldots, q_{p_2+1}')\in Q_{k_2-1,p_2+1},
    \end{aligned}
  \end{equation}
 then we have
\begin{equation}
  h_q(\partial_X)\ h_{q'} (\partial_Y)\ \DS
   = c_{\tilde{\gamma}}\ \Delta_{\tilde{\gamma}} + 
   \sum_{\gamma'' >_{\content} \tilde{\gamma}}
     c_{\gamma''}\ \Delta_{\gamma''}
\end{equation}
where $c_{\tilde{\gamma}} > 0$ and
\begin{equation}\label{eq:maxgam}
  \tilde{\gamma} = \gamma[q_1,q_2,\ldots,q_{p_1+1};q_1',q_2',\ldots,q_{p_2+1}'].
\end{equation}
\end{theorem}

\begin{proof}
  This follows from the argument in the proof of Lemma~\ref{lem:inv}.
  Recall that since a complete symmetric function $h_j$ can be written
  as a nonnegative linear combination of monomial symmetric functions,
  the action of $h_j(\partial_X)$ can be viewed as the movement of
  boxes to the left such that the distances traveled sum to $j$.  In
  the proof of that lemma, we defined an involution that (among other
  things) canceled the contributions coming from any $m_\mu$ with $\mu
  \neq (j)$.

  In the situation of this theorem, we can proceed analogously.  Here,
  though, we need to allow for the possibility that a box moved by an
  $h_{q_i}$ also gets moved by an $h_{q_j}$ for some $j \neq i$.
  However, any such resulting $\gamma''$ will satisfy $\gamma''
  >_{\content} \tilde{\gamma}$.
\end{proof}

The primary goal of this paper is to prove the following theorem.

\begin{theorem}\label{T:bigbasisth}
  The collection $\hhBPgam$ is linearly independent in $\RI$. Hence,
  by the equality of \eqref{E:HilbJ} and \eqref{E:rightsum},
  $\hhBPgam$ forms a basis for $\RI.$
\end{theorem}

By \cite[Theorem 2.1]{Allend}, the subspace of $\BPgam$ spanned by
bipermanents $[T,V]_{\per}$ of a given shape $\lambda$ and given $V$
carries a copy of the irreducible $S^\lambda$.  It
is this fact that lets us conclude Theorem~\ref{T:maintheorem} from
the Hilbert series \eqref{E:HilbJ} of Corollary~\ref{cor:final}.

Note that it is enough to show that the collection
\begin{equation}
  \left\{
  b(\partial_X,\partial_Y) \ \DS: b\in\hhBPgam\right\}
\end{equation}
is linearly independent in $\CXY$ since
\begin{equation}\label{eq:formainthm}
\begin{aligned}
  \sum_{j=1}^i c_j b_j \equiv 0 \pmod{\IS} &\Longleftrightarrow 
  \left(\sum_{j=1}^i c_j b_j \right)(\partial_X,\partial_Y)\ \DS = 0\\
  &\Longleftrightarrow \sum_{j=1}^i 
  \left(c_j\ b_j(\partial_X,\partial_Y)\ \DS\right) = 0.
\end{aligned}
\end{equation}

To do this, however, we need the following lemma.  A proof can be
found in \cite{Allena} (see equations (6.5) and (6.6) in Theorem 6.2),
however, there is a sign missing, so we include a proof here.  The
lemma lets us expand the action of a bipermanent on a determinant
associated to a lattice diagram.  To state the lemma succinctly, we
first associate to the bipermanent in question a family of $\A$-filled
diagrams.  So, consider a bitableau $(T,C)\in\bitabsetn$ and set $S =
\std(C)$.  Write $\iota = \iota_{T,S}$ for the map that takes $i$ to
$s^T_i$.  For any permutation $\phi\in S_n$, we then define an
$\A'$-filled diagram $E_\phi^{\balpha}$ by placing
$\alpha_{\phi(\iota(i))}-c^T_i$ in the cell containing $i$ in
$T$.

\begin{lemma}\label{lem:genderivs}
  Let $\balpha=[\alpha_1,\alpha_2, \ldots, \alpha_n]$ and $(T,C)\in
  \bitabsetn$.  Then
  \begin{equation}\label{eq:tcperlem}
 [T,C]_{\per}(\prt{X}, \prt{Y} )\ \Dalph  = 
    \sgn(\iota) \sum_{\phi\in S_n} \sgn(\phi)\ d_\phi\ [T^t,(E_\phi^{\balpha})^t]_{\det},
  \end{equation}
  for the $\A'$-filled diagrams $E_\phi^{\balpha}$ defined above and integers
  $d_\phi \geq 0$.  We make the convention that $d_\phi = 0$ if any
  entry of $E_\phi^{\balpha} \in \A'-\A$ or if $[T^t,(E_\phi^{\balpha})^t]_{\det} = 0$;
  otherwise $d_\phi > 0$.  
\end{lemma}
Note that $[T^t,(E_\phi^{\balpha})^t]_{\det} = 0$ when there is a repetition in
some row of $E_\phi^{\balpha}$.
\begin{proof}
  Define $q_{T,C} = q_{T,C}(\partial_X,\partial_Y)$ to be the monomial
  $\partial_{z_1}^{c^T_1}\partial_{z_2}^{c^T_2}\cdots
  \partial_{z_n}^{c^T_n}$.  Then $[T,C]_{\per}$ can be
  expanded as $\sum_{\sigma\in R_T} \sigma(q_{T,C})$.  Expanding
  $\Dalph$ as well, we have the following expression for the left-hand
  side of \eqref{eq:tcperlem}:
  \begin{equation}\label{eq:tcpereq}
  \begin{aligned}
    {[T,C]_{\per}}(\partial_X,\partial_Y)\Dalph &= 
    \sum_{\sigma\in R_T}\sum_{\tau\in S_n} \sgn(\tau)
    (\sigma q_{T,C})
    \prod_{j=1}^n z_j^{\alpha_{\tau^{-1}(j)}}\\
    &= \sum_{\sigma\in R_T}\sum_{\tau\in S_n} \sgn(\tau)
    \,\sigma\negmedspace \left(q_{T,C}
    \prod_{j=1}^n z_j^{\alpha_{\sigma(\tau^{-1}(j))}}\right)\\
    &= \sum_{\phi\in S_n} \sum_{\sigma\in R_T} \sgn(\iota\circ\phi\circ\sigma)
    \,\sigma\negmedspace \left(q_{T,C}
    \prod_{j=1}^n z_j^{\alpha_{\phi(\iota(j))}}\right).
  \end{aligned}
  \end{equation}
  As an illustration of the equality between the first two lines, consider
  \begin{equation}
    ((1,5,3)\partial_{x_5}^4) x_3^{\alpha_3} = 
    \alpha_3^{\underline{4}} x_3^{\alpha_3 - 4} 
    = (1,5,3)(\partial_{x_5}^4 x_5^{\alpha_3}).
  \end{equation}
  In going from the second to the third lines, we set $\phi =
  \sigma\circ \tau^{-1}\circ\iota^{-1}$, used the fact that the sign of any
  permutation is the sign of its inverse, and noted that as $\tau$
  runs over $S_n$, so does $\phi$.

  For each $\phi$, we would like to interpret the sum over $\sigma$ as
  a bideterminant $[T^t,(E_\phi^{\balpha})^t]_{\det}$.  Formally this
  makes sense as bideterminants are signed sums over the elements in
  the column stabilizer of some filled diagram.  Here we have a signed
  sum over a row stabilizer; hence we consider transposes.  If
  $q_{T,C}(\partial_X,\partial_Y)$ were not in the equation, we would
  use $\alpha_{\phi(\iota(i))}$ as our entry in $E_\phi^{\balpha}$
  corresponding to $i$ in $T$.  (The $\iota$ accounts for the fact
  that we are not assuming $S = T$.)  Up to the multiplicative
  constants $d_\phi$, the action of $q_{T,C}(\partial_X,\partial_Y)$
  is to subtract $c^T_i$.  This is consistent with the statement of
  the lemma.
\end{proof}

\begin{example}
  We illustrate Lemma~\ref{lem:genderivs} with the following simple
  computation.  In this example, the identity is the only $\phi$ for
  which the entries of $E_\phi^{\balpha}$ are all in $\A$.  Also,
  $\iota$ is given in cycle notation by $(2,4,3)$ which yields
  $\sgn(\iota) = 1$.  In the second-to-last line, the subtraction of
  filled diagrams should be interpreted entrywise.
  \begin{equation}
    \begin{aligned}
    \left[\,\young(4,3,12)\,\right.&\left.,
      \young(\ulb,\ula,\olb\ulc)\ \right] _{\per}(\partial_X,\partial_Y)
    \Dgenalph{(\ol{2},\ul{0},\ul{1},\ul{3})}\\
    &= \left(\partial_{y_1}\partial_{x_2}^2\partial_{x_4} + 
    \partial_{y_2}\partial_{x_1}^2\partial_{x_4}\right)
    \det
    \begin{pmatrix}
      y_1^2 & y_2^2 & y_3^2 & y_4^2\\
      1 & 1 & 1 & 1\\
      x_1 & x_2 & x_3 & x_4\\
      x_1^3 & x_2^3 & x_3^3 & x_4^3
    \end{pmatrix}\\
    &= 12y_1x_2 - 12y_2x_1 = 
    12\,\left[\,\young(2,134)\,,\young(\uld,\olc\ula\ulb) - 
      \young(\ulc,\olb\ula\ulb)\ \right]_{\det}\\
    &= 12\,\left[\,\young(2,134)\,,\young(\ulb,\olb\ula\ula)\ \right]_{\det}.
    \end{aligned}
  \end{equation}
\end{example}

For $U$ an $\A$-filled diagram, we define 
\begin{equation}
  \begin{aligned}
  \conx(U) &= [u_{m_2+p_2+2}, u_{m_2+p_2+3},\ldots,u_{n}],\\
  \cony(U) &= [u_1, u_2, \ldots u_{m_2+p_2}],
  \end{aligned}
\end{equation}
where the square brackets indicate that we have arranged the elements
of the sets in increasing order.  Let $(T_1,U_1), (T_2,U_2)$ be
bitableaux.  We define
\begin{equation}\label{eq:rlexpair}
  (T_1,U_1)\bitabord (T_2,U_2)
\end{equation}
according to the following tiebreakers (recall that $\colseq(U)$
denotes the column sequence of $U$; cf. Example~\ref{ex:rkc}):
\begin{enumerate}
\item $\shape(T_1^t) <_{\normlexord} \shape(T_2^t)$;
\item $\conx(U_1) >_{\normlexord} \conx(U_2)$;
\item $\cony(U_1) <_{\normlexord} \cony(U_2)$;
\item $\colseq(U_1) <_{\lexord{\A}} \colseq(U_2)$;
\item $\colseq(T_1) >_{\lexord{\A}} \colseq(T_2)$. 
\end{enumerate}

We are now ready to prove that $\hhBPgam$ is a basis for $\RI$.  As
outlined in \eqref{eq:formainthm}, it suffices to prove the following
theorem.

\begin{theorem}\label{T:triangular}
  Let $T,U\in\stdtabn$, $q = (q_1, q_2, \ldots, q_{p_1+1})\in
  Q_{k_1-1,p_1+1}$ and $q' = (q_1', q_2', \ldots, q_{p_2+1}')\in
  Q_{k_2-1,p_2+1}$.  Set $C = \cochg{\gamma}(U)$.  Then
  \begin{equation}\label{eq:finthm}
    h_q(\partial_X)\ h_{q'}(\partial_Y)\ [T,C]_{\per}
    (\partial_X,\partial_Y) \ \DSnoxy = \\ 
    d\ [T^t,P^t]_{\det}  + 
    \sum_{(S,W)} d_{S,W} [S^t,W^t]_{\det},
  \end{equation}
  where $d \neq 0$, $P = E_\varepsilon^{\bbeta}$, and the sum is over
  all $(S,W)\in\bitabsetn$ with $(S,W) \bitabordg (T,P)$.
\end{theorem}

Before presenting the proof, we illustrate \eqref{eq:finthm} with an
explicit computation.
\begin{example}\label{ex:threecases}
  Let $\gamma$ derive from the collection of boxes $\balpha =
  (\ol{4},\ol{3},\ul{0},\ul{3},\ul{4},\ul{5})$.  Set $T =
  \young(46,35,12)$ and $U = \young(56,24,13)\,$.  Then $C =
  \cochg{\gamma}(U) = \young(\ulc\ulc,\ula\ulb,\olb\ula)\,$ and
  $\iota = (2,3)(4,5)$.  We consider the expansion of
  \begin{equation}\label{eq:bigexexp}
  \begin{aligned}
    h_{1,1}(\partial_X)& h_{2,1}(\partial_Y)
    [T,C]_{\per}(\partial_X,\partial_Y) \DSnoxy\\ 
    &= [T,C]_{\per}(\partial_X,\partial_Y) 
    h_{1,1}(\partial_X)h_{2,1}(\partial_Y)\DSnoxy\\
    &= [T,C]_{\per}(\partial_X,\partial_Y)
      \left(c_{\bbeta}\ \Dgenalph{\bbeta} + 
      \sum_{\bdelta >_{\content} \bbeta} c_{\bdelta}\ \Dgenalph{\bdelta} \right)\\
    &= c_{\bbeta}\ [T,C]_{\per} (\partial_X,\partial_Y) \ \Dgenalph{\bbeta}  + 
    \sum_{\bdelta >_{\content} \bbeta} c_{\bdelta}\ 
    [T,C]_{\per} (\partial_X,\partial_Y)
    \ \Dgenalph{\bdelta},
  \end{aligned}
  \end{equation}
  where $\bbeta = (\ol{3},\ol{1},\ul{0},\ul{2},\ul{3},\ul{5})$.  For
  $\gamma$ associated to the $\balpha$ above, $m_2 + p_2 + 1 = 3$.  In
  addition
  \begin{equation}
    [T,C]_{\per}(X_6,Y_6) = 
    \sum_{\sigma\in S_{\{1,2\}}\times S_{\{3,5\}} \times S_{\{4,6\}}} 
    z_{\sigma(1)}^{\ol{1}}z_{\sigma(2)}^{\ul{0}}z_{\sigma(3)}^{\ul{0}}
    z_{\sigma(4)}^{\ul{2}}z_{\sigma(5)}^{\ul{1}}z_{\sigma(6)}^{\ul{2}}.
  \end{equation}
  It follows then that
  \begin{equation}
    [T,C]_{\per}(\partial_{X_6},\partial_{Y_6}) = 
    (\partial_{y_1} + \partial_{y_2})(\partial_{x_3} + \partial_{x_5})
    (2\cdot \partial_{x_4}^2\partial_{x_6}^2). 
  \end{equation}
  So, by Lemma~\ref{lem:genderivs},
  $[T,C]_{\per}(\partial_{X_6},\partial_{Y_6}) 
  \Dgenalph{\bbeta}$ can be expanded as 
  \begin{equation}
    \sgn(\iota_{T,U}) \sum_{\phi\in S_6} 
    \sgn(\phi)\, d_\phi [T^t,(E_\phi^{\bbeta})^t]_{\det}.
  \end{equation}

  It is easily checked for this particular example that the entries of
  $E_\phi^{\bbeta}$ will all be in $\A$ only if $\phi \in
  S_{1,2,3}\times S_{4,5,6}$ with $\phi(1) \neq 3$.  There are
  $24$ such permutations.  Table~\ref{tab:ephi} gives these $\phi$
  along with the coefficients $d_\phi$ and diagrams $E_\phi^{\bbeta}$.
  We leave the expansion of the sum in \eqref{eq:bigexexp} to the reader.
  This completes our example.

\newcommand{\minitab}[2][l]{\begin{tabular}{#1}#2\end{tabular}}
\begin{table}[h]
  \centering
  \setlength{\extrarowheight}{10pt}
  \begin{tabular}{|c|>{$}c<{$}>{$}c<{$}>{$}c<{$}!{\vrule width 1pt}
      >{$}c<{$}>{$}c<{$}>{$}c<{$}|}\hline
     & \phi & d_\phi & (E_\phi^{\bbeta})^t & 
     \phi & d_\phi & (E_\phi^{\bbeta})^t \\\hline
     Cases & \varepsilon & 720 & 
     \minitab[c]{\young(\ula\ulb\uld,\olc\olb\ulb)} &
     (5,6) & 720 & \minitab[c]{\young(\ula\ulb\ulb,\olc\olb\uld)} \\ 
     1 \& 2 & (2,3) & 720 & \minitab[c]{\young(\olb\ulb\uld,\olc\ula\ulb)} &
     (2,3)(5,6) & 720 & \minitab[c]{\young(\olb\ulb\ulb,\olc\ula\uld)} \\\hline 
     & (4,5) & 360 & \minitab[c]{\young(\ula\ulc\uld,\olc\olb\ula)} & 
     (4,6) & 120 & \minitab[c]{\young(\ula\ule\ula,\olc\olb\ulb)}\\
     & (2,3)(4,5) & 360 & \minitab[c]{\young(\olb\ulc\uld,\olc\ula\ula)} & 
     (2,3)(4,6) & 120 & \minitab[c]{\young(\olb\ule\ula,\olc\ula\ulb)}\\
     & (4,6,5) & 180 & \minitab[c]{\young(\ula\ule\ulb,\olc\olb\ula)} & 
     (4,5,6) & 360 & \minitab[c]{\young(\ula\ulc\ula,\olc\olb\uld)} \\ 
     & (2,3)(4,6,5) & 180 & \minitab[c]{\young(\olb\ule\ulb,\olc\ula\ula)} & 
     (2,3)(4,5,6) & 360 & \minitab[c]{\young(\olb\ulc\ula,\olc\ula\uld)} \\ 
     Case 3 & (1,2,3) & 240 & \minitab[c]{\young(\old\ulb\uld,\ula\ula\ulb)} & 
     (1,2,3)(4,5) & 120 & \minitab[c]{\young(\old\ulc\uld,\ula\ula\ula)}\\ 
     & (1,2,3)(4,6) & 60 & \minitab[c]{\young(\old\ule\ula,\ula\ula\ulb)} & 
     (1,2,3)(5,6) & 240 & \minitab[c]{\young(\old\ulb\ulb,\ula\ula\uld)}\\ 
     & (1,2,3)(4,6,5) & 60 & \minitab[c]{\young(\old\ule\ulb,\ula\ula\ula)} & 
     (1,2,3)(4,5,6) & 120 & \minitab[c]{\young(\old\ulc\ula,\ula\ula\uld)}\\ 
     & (1,2) & 0 & \minitab[c]{\young(\ula\ulb\uld,\ula\old\ulb)} & 
     (1,2)(4,5) & 0 & \minitab[c]{\young(\ula\ulc\uld,\ula\old\ula)}\\ 
     & (1,2)(4,6) & 0 & \minitab[c]{\young(\ula\ule\ula,\ula\old\ulb)} & 
     (1,2)(5,6) & 0 & \minitab[c]{\young(\ula\ulb\ulb,\ula\old\uld)}\\ 
     & (1,2)(4,6,5) & 0 & \minitab[c]{\young(\ula\ule\ulb,\ula\old\ula)} & 
     (1,2)(4,5,6) & 0 & \minitab[c]{\young(\ula\ulc\ula,\ula\old\uld)}\\\hline 
  \end{tabular}
  \caption{Expansion of Example~\ref{ex:threecases}}
  \label{tab:ephi}
\end{table}
\end{example}

\begin{proof}
  Let $T,U\in\stdtabn$.  Using Theorem~\ref{thm:hdelta}, we have that
  \begin{equation}
    \begin{aligned}
    h_q(\partial_X)&h_{q'}(\partial_Y)\ [T,C]_{\per}
    (\partial_X,\partial_Y) \ \DS\\  
    &= [T,C]_{\per} (\partial_X,\partial_Y) 
      h_q(\partial_X)h_{q'}(\partial_Y)\ \DS\\
    &= [T,C]_{\per} (\partial_X,\partial_Y) 
      \left(c_{\bbeta}\ \Dgenalph{\bbeta} + 
      \sum_{\bdelta >_{\content} \bbeta} c_{\bdelta}\ \Dgenalph{\bdelta} \right)\\
    &= c_{\bbeta}\ [T,C]_{\per} (\partial_X,\partial_Y) \ \Dgenalph{\bbeta}  + 
    \sum_{\bdelta >_{\content} \bbeta} c_{\bdelta}\ [T,C]_{\per} (\partial_X,\partial_Y)
    \ \Dgenalph{\bdelta},
    \end{aligned}
  \end{equation}
  with $c_{\bbeta} > 0$.  Now (abbreviating $\iota_{T,U}$ by
  $\iota$), Lemma~\ref{lem:genderivs} implies that
  \begin{equation} [T,C]_{\per}( \prt{X} , \prt{Y} )\ \Dgenalph{\bbeta}  = 
    \sgn(\iota) \sum_{\phi\in S_n} \sgn(\phi)\ d_\phi\ 
    [T^t,(E_\phi^{\bbeta})^t]_{\det}.
  \end{equation}
  Therefore,
  \begin{multline}\label{E:Ephidelta}
    h_q(\partial_X)\,h_{q'}(\partial_Y)[T,C]_{\per}
    (\partial_X,\partial_Y) \ \DS\\  
    = \sgn(\iota) \sum_{\phi\in S_n}
      \sgn(\phi)\ \left( c_{\bbeta}\ d_\phi\ [T^t,(E_\phi^{\bbeta})^t]_{\det}  + 
        \sum_{\delta >_{\content} \beta} c_{\bdelta}\
        d_\phi\ [T^t,(E_\phi^{\bdelta})^t]_{\det}\right).
  \end{multline}
  Our goal is to show that there is a minimum bitableau (with respect
  to $\bitabord$) occurring on the right-hand side of
  \eqref{E:Ephidelta}.  First note that many of the bitableaux
  $(T^t,(E_\phi^{\bbeta})^t)$ and $(T^t,(E_\phi^{\bdelta})^t)$ are not
  standard.  However, Theorem~\ref{T:Rotaexp} tells us how the shapes
  and column sequences are affected by straightening.

  We split into three cases dependent on the indexing permutations
  $\phi$.  (A proof of an argument with similar statements and
  complete details can be in found in \cite[Theorem 6.2]{Allena}.)
  \begin{enumerate}
  \item[Case 1:] $\phi = \varepsilon$.\\
    We get the first term on the right-hand side of \eqref{eq:finthm}
    with $d = c_{\bbeta}\, d_\varepsilon > 0$.

  \item[Case 2:] $\phi\neq \varepsilon$, but $\phi(z_1^{c^T_1}\cdots
    z_n^{c^T_n}) =
    z_1^{c^T_1}\cdots z_n^{c^T_n}$.\\
    For such $\phi$, $\kappa(E_\phi^{\bbeta}) = \kappa(P)$.  However,
    it follows from the construction of $C = \cochg{\gamma}(U)$ from
    $U$ that $\rowseq((E_\phi^{\bbeta})^t) = \colseq(E_\phi^{\bbeta})
    >_{\lexord{\A}} \colseq(P) = \rowseq(P^t)$.  Note that the
    $E_\phi^{\bbeta}$ arising in this case are standard and do not
    need to be straightened.

  \item[Case 3:] $\phi(z_1^{c^T_1}\cdots z_n^{c^T_n}) \neq 
    z_1^{c^T_1}\cdots z_n^{c^T_n}$.\\
    We must have $\conx(E_\phi^{\bbeta}) <_{\normlexord} \conx(P)$ or
    $\cony(E_\phi^{\bbeta}) >_{\normlexord} \cony(P)$.
  \end{enumerate}

  It is not difficult to see that the $E_\phi^{\bdelta}$ with $\bdelta
  >_{\content} \bbeta$ that
  \begin{equation}
    (T,E_\phi^{\bdelta}) \bitabordg (T,P).
  \end{equation}
  Thus, combinations of $\phi$ and $\bdelta\neq \bbeta$ substituted
  into equation \eqref{E:Ephidelta} (and using Theorem
  \ref{T:Rotaexp}) yield bideterminants $[S,W]_{\det}$ with $(S,W)
  \bitabordg (T,P)$.  Thus, we have
  \begin{multline}
    h_q(\partial_X)h_{q'}(\partial_Y)\ [T,C]_{\per}
    (\partial_X,\partial_Y) \ \DS = \\ 
    d\ [T^t,P^t]_{\det}  + 
    \sum_{(S,W) \bitabordg (T,P)} d_{S,W} [S^t,W^t]_{\det},
  \end{multline}
  with $d \neq 0$.  This proves the theorem.
\end{proof}

Since we have that $\JS \subset \IS$, the fact that the Hilbert series
in equation \eqref{E:HilbJ} equals the summation in equation
\eqref{E:rightsum} implies that we must have $\JS=\IS$. Thus, equation
\eqref{E:HilbJ} must give the Hilbert series for $\RI$. Since the
collection $\hhBPgam$ is linearly independent in $\RI$ and it gives the
correct Hilbert series, we must have that $\hhBPgam$ is in fact a basis for
$\RI$.  The proof of Theorem \ref{T:bigbasisth} is now complete.

\section{Some notes, applications and conjectures}

\begin{note}
  Suppose $\mathcal{D}$ is any basis for $\RI$ where $\gamma =
  ((m_1,n-m_1+1),(0,0),(0,0))$ for some $m_1 > 0$.  The basis
  $\mathcal{D}$ can be substituted in the place of $\BPgam$ in the
  definition of $\hhBPgam$ of equation \eqref{D:basisdef} such that
  $\hhBPgam$ still yields a basis for $\RI$.  Examples of such bases
  $\mathcal{D}$ include descent monomials (see \cite{Allena},
  \cite{Garsia} or \cite{Reiner}), Artin monomials (see \cite{Artin}),
  Schubert monomials (see \cite{Schubert}) and Higher Specht
  Polynomials (see \cite{Morita}).
\end{note}

\begin{note}
  The ideas of the previous sections are
  easily extendable to the complex reflection groups $G(r,p,n)$.  In
  this case, the lattice diagrams utilized are those of the form
  $\latdiag{\bbeta}$ where $\latdiag{\balpha}$ is a hollow lattice
  diagram and $\bbeta = \{(r\alpha_{i,1},r\alpha_{i,2}): \alpha\in
  \balpha\}$.  The resulting bases give representations of the complex
  reflection groups by ways of $m$-tableaux, standard $m$-tableaux and
  cocharge $m$-tableaux.  See \cite{Allena} to see how this
  translation is accomplished.
\end{note}

\begin{note}
  Some of the results of this paper can be extended to diagonally
  symmetric and anti-symmetric rings in the four sets of variables
  $X_n$, $Y_n$, $Z_n$ and $W_n$.  The diagonal action of $S_n$ on
  $\bbC[X_n,Y_n,Z_n,W_n]$ and the rings of symmetric polynomials and
  anti-symmetric polynomials $\bbC^+[X_n,Y_n,Z_n,W_n]$ and
  $\bbC^-[X_n,Y_n,Z_n,W_n]$ are defined in the natural manner.  For
  $\sigma\in S_n$, set $\varepsilon^+(\sigma) = 1$,
  $\varepsilon^-(\sigma) = \sgn(\sigma)$.  For each of $+$ and $-$,
  define
  \begin{multline}
    R_{\gamma_1,\gamma_2}^{\pm} = 
    \C^{\pm}[X_n,Y_n,Z_n,W_n]/\bigl\{P\in \bbC^{\pm}[X_n,Y_n,Z_n,W_n]:\\
      P(\partial_X,\partial_Y,\partial_Z,\partial_W) \
      \Delta_{\gamma_1}(X_n,Y_n)\ \Delta_{\gamma_2}(Z_n,W_n) = 0\bigr\}.
  \end{multline}
  Furthermore, for an arbitrary standard tableau $Q$, set
    \begin{align}
    [T_1,T_2]_{\per}^+ &= \sum_{\sigma\in S_n} \sigma
    \left([Q,T_1]_{\per}(X_n,Y_n)\ [Q,T_2]_{\per}(Z_n,W_n)\right)\\
    [T_1,T_2]_{\per}^- &= \sum_{\sigma\in S_n} \sigma
    \left([Q,T_1]_{\per}(X_n,Y_n)\ [Q,T_2]_{\per}(\partial_Z,\partial_W)
      \Delta_{\gamma_2}(Z_n,W_n)\right).  
    \end{align}
  Using the techniques found in \cite{Allenb}, it is not difficult to
  show the following theorem:
  \begin{theorem}
    For each of $+$ and $-$,
    \begin{multline}
      \hhBPgamii^{\pm} = 
      \left\{ h_{q_1, q_2, \ldots, q_{p_1+1}}(X,Y)\ 
        h_{q_1', q_2', \ldots, q_{p_1+1}'}(Z,W)\ [T_1,T_2]_{per}^{\pm}:\right.\\
      \left.\shape(T_1) = \shape(T_2), T_1 \in \cochgalln{\gamma_1}, 
        T_2\in \cochgalln{\gamma_2}\right\},
    \end{multline}
    where $(q_1, q_2, \ldots, q_{p_1+1})\in Q_{k_1-1,p_1+1}$ and
    $(q_1', q_2', \ldots, q_{p_2+1}')\in Q_{k_2-1,p_1+1}$, is a basis
    for $R_{\gamma_1,\gamma_2}^{\pm}$.
  \end{theorem}

  These are generalizations of rings that have been studied in
  \cite{Allenc,Reiner}, for example.  Furthermore, it should be noted
  that this theorem is not a generalization of the theorems found in
  \cite{Allenb}. Rather, we are just noting that the techniques found
  there easily generalize to the situation in this paper.
\end{note}

\bigskip
\bibliographystyle{acm}
\begin{small}
\bibliography{refs}
\end{small}

\end{document}

%% file: youngdiag.tex
   \newdimen\Squaresize \Squaresize=10pt
   \newdimen\Thickness \Thickness=0.5pt

   \def\Square#1{\hbox{\vrule width \Thickness
        \vbox to \Squaresize{\hrule height \Thickness\vss
           \hbox to \Squaresize{\hss#1\hss}
        \vss\hrule height\Thickness}
   \unskip\vrule width \Thickness}
   \kern-\Thickness}

   \def\Vsquare#1{\vbox{\Square{$#1$}}\kern-\Thickness}

   \newdimen\squaresize \squaresize=30pt
   \newdimen\thickness \thickness=0.2pt

   \def\SSquare#1{\hbox{\vrule width \thickness
        \vbox to \squaresize{\hrule height \thickness\vss
           \hbox to \squaresize{\hss#1\hss}
        \vss\hrule height\thickness}
   \unskip\vrule width \thickness}
   \kern-\thickness}

   \def\vsquare#1{\vbox{\SSquare{$#1$}}\kern-\thickness}

   \def\thisbox#1{\kern-.09ex\fbox{#1}}
   \def\downbox#1{\lower1.200em\hbox{#1}}
